\documentclass[11pt]{article}
\usepackage{graphicx}
\usepackage{epstopdf}
\usepackage{wrapfig}

\usepackage{amsmath,amsfonts,amssymb}
\def\tto{\;{\lower 1pt \hbox{$\rightarrow$}}\kern -10pt
\hbox{\raise 2pt \hbox{$\rightarrow$}}\;}

\def\Tilde{\widetilde}

\def\ra{\rangle}
\def\la{\langle}

\def\B{I\!\!B}
\def\h{\hfill\square}
\def\R{I\!\!R}

\def\ox{\bar{x}}

\newcommand{\chuan}[1]{\rVert #1 \rVert}

\def\h{\hfill\square}

\def \R{\Bbb R}

\newcounter{lk}

\begin{document}
\begin{center}
\vspace*{0.3in} {\bf SOLUTIONS CONSTRUCTIONS OF A GENERALIZED SYLVESTER PROBLEM AND A GENERALIZED FERMAT-TORRICELLI PROBLEM FOR EUCLIDEAN BALLS}\\[2ex]
Nguyen Mau Nam,\footnote{Fariborz Maseeh Department of Mathematics and Statistics,
Portland State University, Portland, OR 97202, United States, email:
mau.nam.nguyen@pdx.edu (corresponding author).}
Nguyen Hoang,\footnote{Department of Mathematics, College of
Education, Hue University, Hue City, Vietnam, email:
nguyenhoanghue@gmail.com.} and Nguyen Thai An\footnote{Institute of Mathematics, Vietnam Academy of Science and Technology, 18 Hoang Quoc Viet,
Hanoi 10307, Vietnam, email: thaian2784@gmail.com. }

\end{center}
\small{\bf Abstract:}  The classical Apollonius' problem is to
construct circles that are tangent to three given circles in a
plane. This problem was posed by Apollonius of Perga in his work ``Tangencies". The Sylvester problem, which
was introduced by the English mathematician J.J. Sylvester, asks for the smallest circle that encloses a finite
collection of points in the plane. In this paper, we study the
following generalized version of the Sylvester problem and its
connection to the problem of Apollonius: given two finite
collections of Euclidean balls in $\Bbb R^n$, find the smallest Euclidean
ball that encloses all of the balls in the first collection and
intersects all of the balls in the second collection. We also study a generalized version of the Fermat-Torricelli problem stated as follows: given two finite
collections composed of three Euclidean balls in $\Bbb R^n$, find a point that minimizes the sum of the farthest distances to the balls in the first collection and shortest distances to the balls in the second collection.

\medskip
\vspace*{0,05in} \noindent {\bf Key words.}  Convex analysis and
optimization, generalized differentiation, smallest enclosing circle
problem, Fermat-Torricelli problem.

\newtheorem{theorem}{Theorem}[section]
\newtheorem{proposition}{Proposition}[section]
\newtheorem{lemma}{Lemma}[section]
\newtheorem{corollary}{Corollary}[section]
\newtheorem{definition}{Definition}[section]
\newtheorem{example}{Example}[section]
\newtheorem{remark}{Remark}[section]

\renewcommand{\theequation}{\thesection.\arabic{equation}}
\normalsize

\section{Introduction}

The celebrated \emph{Sylvester problem} asks for the smallest circle that encloses a finite number of given points in a plane. This problem was introduced by the English mathematician James Joseph Sylvester in \cite{syl}. Because of its importance in many different applications, the problem has attracted many researchers from different fields. The Sylvester problem and its version in higher dimensions are now called under different names such as: the \emph{smallest enclosing ball problem}, the \emph{minimum ball problem}, or the \emph{bomb problem}. The readers are referred to [2-7] and the references therein for recent study on the problem and its generalizations.\vspace*{0.05in}

In the 17th century, the French mathematician Pierre de Fermat introduced another problem of Euclidean geometry asking for a point that minimizes the sum of the distances to three given points in the plane. The problem was then solved by the Italian mathematician and physicist Evangelista Torricelli, and it is now known as the \emph{Fermat-Torricelli problem}. Numerous articles have been written to study this problem and its generalizations to many different settings; see, e.g., [8-12].\vspace*{0.05in}

In our recent publications [13-16], we study extensions of the Sylvester and Fermat-Torricelli problems in which the given points are replaced by given sets. The existence and uniqueness of optimal solutions, optimality conditions, numerical algorithms, as well as other properties of the problems have been initially addressed. In particular, in \cite{nh}, we study of the following problem called the \emph{generalized Sylvester problem}: given two finite
collections of sets in a normed space, find a ball of the smallest radius, whose center lies
in a given constraint set, that encloses all
the sets in the first collection and intersects all the sets in the
second one. We also introduce and study the following generalized version
of the classical Fermat-Torricelli problem called the
\emph{generalized Fermat-Torricelli problem}: given two finite
collections of sets in a normed space, find a point in a given
constraint set that minimizes the sum of the \emph{farthest distances} to the sets in the first collection and \emph{shortest
distances/distances} to the sets in the second
collection.\vspace*{0.05in}

This paper is a continuation of our development with the study of a special case of the unconstrained generalized Sylvester and Fermat-Torricelli problems in which the given sets are Euclidean balls. The special features of the Euclidean balls therein make the problems distinct from the general case. Although the new problems seem to be very interesting, to the best of our knowledge, they have not been considered in the literature. We will develop an approach in which the tools of modern convex analysis and optimization are employed to solve the problem. Our paper is organized as follows. In section 2, we focus on theoretical study of the generalized Sylvester problem for Euclidean balls. An important specification for the case of three balls and its connection to the problem of Apollonius are carefully investigated. The construction of solutions is presented for each specific problem. In section 3, we study the generalized Fermat-Torricelli problem for three Euclidean balls. We establish a necessary and sufficient condition for the existence and uniqueness solution of the problem in specific cases. The construction of the solutions is also given.

\section{Preliminaries}
\label{s:Pre}

In this section, we will formulate mathematical models for the problems under consideration and present some concepts and results of convex analysis that will be used in the next sections.

Let $I$ and $J$ be two finite index sets such that $|I|+|J|>1$. Let $\B$ be the unit ball and let $\Omega_i:=\B(a_i; r_i)$ ($r_i\geq 0$) for $i\in I$ and $\Theta_j:=\B(b_j; s_j)$ ($s_j\geq 0)$ for $j\in J$
be two collections of closed balls in $\R^n$ with the Euclidean norm $\|\cdot\|$. These are our \emph{standing assumptions} throughout the paper unless otherwise stated.

For a closed, bounded and convex set $Q$, the \emph{farthest distance function} and the \emph{distance function} to $Q$ are given respectively by
\begin{equation*}
M(x;Q):=\sup\{\|x-q\|\; |\; q\in Q\}
\end{equation*}
and
\begin{equation*}
D(x;Q):=\inf\{\|x-q\|\; |\; q\in Q\}.
\end{equation*}
The generalized Sylvester problem for these Euclidean balls can be reduced to the following optimization problem:
\begin{equation}\label{modeling}
\mbox{minimize }\mathcal{G}(x):=\max\{ \mathcal{M}(x), \mathcal{D}(x)\}, \; x\in \R^n,
\end{equation}
where
\begin{equation*}
\mathcal{M}(x):=\max\{M(x;\Omega_i)\; |\; i\in I\} \mbox{ and
}\mathcal{D}(x):=\max\{D(x;\Theta_j)\; |\; j\in J\}.
\end{equation*}
Similarly, the mathematical optimization modeling of the generalized Fermat-Torricelli problem is
\begin{equation}\label{F-T0}
\mbox{minimize }\mathcal{H}(x):=\mathcal{M}_1(x)+ \mathcal{D}_1(x), \;x\in \R^n,
\end{equation}
where
\begin{equation*}
\mathcal{M}_1(x):=\sum_{i\in I}M(x;\Omega_i) \mbox{ and
}\mathcal{D}_1(x):=\sum_{j\in J}D(x;\Theta_j).
\end{equation*}

Let us recall in what follows some important concepts and results from convex analysis that will be used throughout the paper. The readers are referred to the books \cite{HU,r} for more systematic development of the field.

A function $\psi: \Bbb R^n\to \R$ is called \emph{convex} iff for any $x, y\in \Bbb R^n$ and for any $\alpha\in ]0,1[$, one has
\begin{equation*}
\psi(\alpha x+(1-\alpha)y)\leq \alpha \psi(x)+(1-\alpha)\psi(y).
\end{equation*}
If this inequality becomes strict whenever $x\neq y$, the function is called \emph{strictly convex}.

A vector $v\in \Bbb \R^n$ is called a \emph{subgradient} of a convex function $\psi$ at $\ox$ iff
\begin{equation*}
\la v, x-\ox\ra \leq \psi(x)-\psi(\ox)\mbox{ for all }x\in \Bbb R^n.
\end{equation*}

The collection of all subgradients of $\psi$ at $\ox$ is called the \emph{subdifferential} of the function at this point and is denoted by $\partial\psi(\ox)$. In the case where $\psi$ is Fr\'echet differentiable at $\ox$, the subdifferential $\partial \psi(\ox)$ reduces to the gradient $\nabla\psi(\ox)$ of the function at $\ox$, and it is a set in the general case.

The following generalization of the classical Fermat's rule called \emph{Fermat subdifferential rule} will be important in the sequel:
\begin{equation}\label{fermat}
\ox \;\mbox{\rm is an absolute minimum of a convex function}\;\psi \;\mbox{\rm if and only if }0\in \partial\psi(\ox).
\end{equation}
Since the functions $\mathcal{G}$ and $\mathcal H$  in problems (\ref{modeling}) and (\ref{F-T0}) are represented respectively in terms of the ``max'' and ``sum" of a finite number of convex functions, we are going to use available subdifferential rules from convex analysis to further explore these problems. If $\psi(x):=\max\{\psi_i(x)\;|\; i=1, \ldots,k\}$, where $\psi_i: \Bbb R^n\to \Bbb R$ for $i=1,\ldots, k$ are convex functions, then
\begin{equation}\label{max}
\partial \psi(\ox)=\mbox{\rm co }\{\partial \psi_i(\ox)\; |\; i\in I(\ox)\},
\end{equation}
where $I(\ox):=\{i=1, \ldots, k\;|\; \psi_i(\ox)=\psi(\ox)\}$ is the active index set at $\ox$. Another important subdifferential rule called the \emph{subdifferential sum rule} for $\psi(x):=\sum_{i=1}^k\psi_i(x)$, is stated as:
\begin{equation}\label{sum}
\partial \psi(\ox)=\sum_{i=1}^k\partial \psi_i(\ox).
\end{equation}

Throughout the paper, we use the following standard notations:
$\mbox{\rm bd }\Omega$ and $\mbox{\rm int }\Omega$ denote
respectively the boundary and the interior of a set $\Omega$; for $a,
b\in \R^n$ and $a\neq b$,
\begin{align*}
&[a, b]:=\{ta+(1-t)b\; |\; t\in [0,1]\};\\
&(a, b):=\{ta+(1-t)b\; |\; t\in\  ]0,1[\};\\
&L(a, b):=\{ta+(1-t)b\; |\; t\in \R\}.
\end{align*}

\section{A Generalized Sylvester Problem for Euclidean Balls}
\label{s:SCNC}

\subsection{Existence and Uniqueness of Optimal Solutions and Optimality Conditions}

Let us start the section with simple formulas for computing distances to Euclidean balls in $\Bbb R^n$ as well as their subdifferentials in the sense of convex analysis; see, e.g., \cite{HU}.

\begin{proposition}\label{subr} Let $\Omega=\B(c;r)$, where $c\in\R^n$ and $r\geq
0$. For any $x\in \R^n$, one has
\begin{equation*}
D(x;\Omega)=\begin{cases}
0, &\text{if }\;\|x-c\|\leq r, \\
\|x-c\|-r, & \text{otherwise},
\end{cases}
\end{equation*}
and $$M(x; \Omega)=\|x-c\|+r.$$
Moreover,
\begin{equation*}
\partial D(\ox;\Omega)=\begin{cases}
N(\ox; \Omega)\cap \B, &\text{if }\;\|\ox-c\|\leq r, \\
\dfrac{\ox-c}{\|\ox-c\|}, & \text{otherwise},
\end{cases}
\end{equation*}
and
\begin{equation*}
\partial M(\ox;\Omega)=\begin{cases}
\B, &\text{if }\;\ox=c, \\
\dfrac{\ox-c}{\|\ox-c\|}, & \text{otherwise}.
\end{cases}
\end{equation*}
\end{proposition}
For any $x\in \Bbb R^n$, define
\begin{equation*}
K(x):=\{i\in I\; |\; M(x;\Omega_i)=\mathcal{G}(x)\} \mbox{ and } L(x):=\{j\in J\; |\; D(x; \Theta_j)=\mathcal{G}(x)\}.
\end{equation*}
The active index set $A(x)$ is given by the \emph{disjoint union} $A(x)=K(x)\cup L(x)$. It is obvious that for every $x\in \Bbb R^n$, one has $1\leq |A(x)|\leq |I|+|J|$.\vspace*{0.05in}

The following theorem establishes a necessary and sufficient condition for the uniqueness of an optimal solution. We will provide a simple direct proof for the result; see also \cite{nh}.
\begin{theorem}\label{eu} The optimization problem {\rm(\ref{modeling})} always has an optimal solution. Moreover, the solution is unique if and only if $I\neq\emptyset$, or $I=\emptyset$ and $\cap_{j\in J}\Theta_j$ contains no more than one point.
\end{theorem}
\noindent {\bf Proof. }The fact that problem {\rm(\ref{modeling})} always has an optimal solution follows from the continuity of the function $\mathcal{G}$ therein and the boundedness of its level sets.

Let us prove the sufficient condition for the uniqueness of an optimal solution. Consider the case where $I\neq \emptyset$. We will first show that
\begin{equation}\label{iq1}
\mathcal{G}(x)=\max\{\|x-a_i\|+r_i, \|x-b_j\|-s_j: i\in I, j\in J\}.
\end{equation}
Since $M(x;\Omega_i)=\|x-a_i\|+r_i$ and $D(x; \Theta_j)=\max\{\|x-b_j\|-s_j, 0\}\geq \|x-b_j\|-s_j$,
one always has
\begin{equation*}
\mathcal{G}(x)\geq\max\{\|x-a_i\|+r_i, \|x-b_j\|-s_j\; |\;  i\in I, j\in J\}.
\end{equation*}
Fix $i_0\in I$. If $K(x)\neq \emptyset$, then for an element $i\in K(x)$, one has
\begin{equation*}
\mathcal{G}(x)=\|x-a_i\|+r_i\leq \max\{\|x-a_i\|+r_i, \|x-b_j\|-s_j\; |\; i\in I, j\in J\}.
\end{equation*}
In the case where $K(x)=\emptyset$, one has $L(x)\neq\emptyset$. Fix $j\in L(x)$. Then $\mathcal{G}(x)=D(x; \Theta_j) > M(x; \Omega_{i_0})\geq 0$. Thus $x\notin \Theta_j$, and hence
\begin{equation*}
\mathcal{G}(x)=\|x-b_j\|-r_j\leq \max\{\|x-a_i\|+r_i, \|x-b_j\|-s_j\; |\; i\in I, j\in J\}.
\end{equation*}
We have justified (\ref{iq1}). Choose a constant $\ell$ such that $\ell >\max\{s_j\; |\; j\in J\}$. It follows from (\ref{iq1}) that $\ox$ is an optimal solution of problem (\ref{modeling}) if and only if, it is a solution of the following optimization problem:
\begin{equation}\label{nmodeling}
\mbox{\rm minimize }\mathcal{\Tilde{G}}(x), \; x\in \Bbb R^n,
\end{equation}
where
\begin{equation*}
\mathcal{\Tilde{G}}(x):=\max\{(\|x-a_i\|+r_i+\ell)^2, (\|x-b_j\|+\ell-s_j)^2\; |\;  i\in I, j\in J\}=(\mathcal{G}(x)+\ell)^2.
\end{equation*}
Since $f_i(x):=(\|x-a_i\|+r_i+\ell)^2$ and $g_j(x):=(\|x-b_j\|+\ell-s_j)^2$ for $i\in I$ and $j\in J$ are strictly convex functions, we see that problem (\ref{nmodeling}) has a unique optimal solution, and hence (\ref{modeling}) also has a unique optimal solution.

Suppose that $I=\emptyset$ and $\cap_{j\in J}\Theta_j$ contains at most one point. If $\cap_{j\in J}\Theta_j$ contains exactly one point $x_0$, then $\mathcal{G}(x_0)=0$ and $x_0$ is the unique solution. In the case $\cap_{j\in J}\Theta_j=\emptyset$  it is also not hard to show that (\ref{iq1}) is satisfied for every $x\in\Bbb R^n$, and problem (\ref{modeling}) has a unique solution.

Let us now prove the necessary condition for the unique of solution. Suppose that problem (\ref{modeling}) has a unique solution and assume by contradiction that $I=\emptyset$ and $\cap_{j\in J}\Theta_j$ contains more than one points. It is obvious that $\mathcal{G}(x)=0$ for any $x\in \cap_{j\in J}\Theta_j$. Thus, any $x\in \cap_{j\in J}\Theta_j$ is an optimal solution of the problem. So we have arrived at a contradiction. The proof is now complete. $\h$\vspace*{0.05in}

For any point $x\in \Bbb R^n$, the \emph{farthest projection} and the \emph{shortest projection} from $x$ to a set $Q$ is given by
\begin{equation*}
\mathcal{P}(x;Q):=\{q\in Q\; |\; \|x-q\|=M(x;Q) \} \mbox{ and }\Pi(x;Q):=\{q\in Q\; |\; \|x-q\|=D(x;Q)\}.
\end{equation*}
If $Q=\B(c; r)$, then
\begin{equation*}
\mathcal{P}(x;Q)=\big\{c-r\dfrac{x-c}{\|x-c\|}\big\},
\end{equation*}
and
\begin{equation*}
\Pi(x;Q)=\begin{cases}
\{x\}, &\text{if }\;\|x-c\|\leq r, \\
c+r\dfrac{x-c}{\|x-c\|}, & \text{otherwise}.
\end{cases}
\end{equation*}
In the case $r_i=0$, the ball $\B(a_i; r_i)=\{a_i\}$, so a ball covers $\B(a_i; r_i)$ if and only if it intersects the ball. Moreover, the case where $I=\emptyset$ has been considered in \cite{naj}. Thus, we exclude these cases in the theorem below for simplicity.

\begin{theorem}\label{optimality} Suppose that $I\neq \emptyset$ with $r_i>0$ for every $i\in I$. An element $\ox$ is the optimal solution of problem {\rm(\ref{modeling})} if and only if one of the following conditions holds:\\[1ex]
{\rm (1)} $\B(\ox; r)$, $r=\mathcal{G}(\ox)$, coincides with $k$ balls from $\Omega_i$ for $i\in I$, $k\geq 1$, $\B(\ox; r)$ contains the other balls in $\{\Omega_i: i\in I\}$, and it intersects $\Theta_j$ for $j\in J$.\\
{\rm (2)} $\mathcal{P}(\ox; \Omega_i)$ and $\Pi(\ox; \Theta_j)$, $i\in K(\ox)$ and $j\in L(\ox)$, are singletons. Moreover, for $p_i:=\mathcal{P}(\ox; \Omega_i)$, $q_i:=\Pi (\ox; \Theta_j)$, $i\in K(\ox)$, $j\in L(\ox)$, one has
\begin{equation*}
\ox\in\mbox{\rm co }\{p_i, q_i\; |\;  i\in K(\ox), j\in L(\ox)\}.
\end{equation*}
\end{theorem}
\noindent {\bf Proof. }Notice that in this case, problem (\ref{modeling}) has a unique solution. By the subdifferential Fermat rule, $\ox$ is the optimal solution of the problem if and only if
\begin{equation}\label{eq2}
0\in \partial \mathcal G(\ox)=\mbox{\rm co }\{\partial M(\ox; \Omega_i), \partial D(\ox; \Theta_j)\; |\; i\in K(\ox), j\in L(\ox)\}.
\end{equation}
Let $\mathcal{J}:=\{i\in K(\ox)\; |\; \partial M(\ox; \Omega_i)\;\mbox{\rm is not a singleton}\}$. If $\mathcal{J}\neq\emptyset$, then $\ox=a_i$  and $\mathcal{G}(\ox)=M(\ox; \Omega_i)=r_i$ for every $i\in \mathcal{J}$. So all balls $\Omega_i$ for $i\in \mathcal{J}$ coincide. Moreover,
\begin{equation*}
r_i=M(\ox; \Omega_i)\geq M(\ox; \Omega_j)=\|a_i-a_j\|+r_j \mbox{ for }i\in \mathcal{J} \mbox{ and }j\in I\setminus \mathcal{J},
\end{equation*}
and $r_i=M(\ox; \Omega_i)\geq D(\ox; \Theta_j)=\|a_i-b_j\|-s_j$ for every $i\in \mathcal{J}$ and $j\in J$. Thus, $\B(\ox; r)$ contains balls in $\{\Omega_i\;|\; i\in I\setminus \mathcal{J}\}$ and intersects $\Theta_j$ for $j\in J$.

Consider the case where $\mathcal{J}=\emptyset$. Then $\partial M(\ox;\Omega_i)$ is a singleton for every $i\in K(\ox)$. This implies $\mathcal P(\ox; \Omega_i)$ is a singleton for such an $i$. For every $j\in L(\ox)$, since $D(\ox; \Theta_j)\geq M(\ox; \Omega_{i_0})>0$ for a fixed index $i_0\in I$, one has that $\ox\notin\Theta_j$, and hence ${\Pi}(\ox;\Theta_j)$ is a singleton. Then (\ref{eq2}) can be equivalently written as
\begin{equation*}
0\in \mbox{\rm co }\{\dfrac{\ox-p_i}{r}, \dfrac{\ox-q_j}{r}\; |\; i\in K(\ox), j\in L(\ox)\}, \ \text {where}\ r=\mathcal{G}(\ox)=\|\ox -p_i\|=\|\ox-q_j\|.
\end{equation*}
Thus, there exist $\lambda_i\geq 0$ for $i\in K(\ox)$ and $\mu_j\geq 0$ for $j\in L(\ox)$ such that $\sum_{i\in K(\ox)}\lambda_i+\sum_{j\in L(\ox)}\mu_j=1$ and
\begin{equation*}
0=\sum_{i\in K(\ox)}\lambda_i \dfrac{\ox-p_i}{r}+ \sum_{j\in L(\ox)}\mu_j \dfrac{\ox-q_j}{r}.
\end{equation*}
Equivalently,
\begin{equation*}
\ox=  \sum_{i\in K(\ox)}\lambda_i p_i+ \sum_{j\in L(\ox)}\mu_j q_j\in \mbox{\rm co }\{p_i, q_i\; |\;  i\in K(\ox), j\in L(\ox)\}.
\end{equation*}
The proof of the converse follows from the Fermat subdifferential rule (\ref{fermat}) since we can verify that (\ref{eq2}) is satisfied in each case. $\h$ \vspace*{0.05in}

Finally, using the well-known \emph{Caratheodory theorem} (see, e.g., \cite[Theorem~1.3.6]{HU}) and Theorem \ref{optimality}, we can prove the proposition below.

\begin{proposition}\label{reduction} Suppose that $I\neq \emptyset$ with $r_i>0$ for every $i\in I$. If $\ox$ is the optimal solution of problem {\rm(\ref{modeling})}, then there exist subindex sets $I_1\subseteq I$ and $J_1\subseteq J$ with $1<|I_1|+|J_1|\leq n+1$ such that $\ox$ is a solution of the generalized Sylvester problem with target sets $\Omega_i$, $i\in I_1$, and $\Theta_j$, $j\in J_1$.
\end{proposition}

\subsection{Three-Ball Generalized Sylvester Problem and the Problem of \\Apollonius}

 In this subsection, we focus on the generalized Sylvester problem for the case of three Euclidean balls in two dimensions. From Proposition \ref{reduction}, we see that this is one of the most important cases since it is possible to reduce the problem with large number of balls to the problem of three balls or less. This observation has been used as the key point for many algorithms to solve the classical minimum ball problem; see, e.g., \cite{wel}.\vspace*{0.05in}

\noindent {\bf Three-Ball Problem: Model I.} The first model we study in this subsection is stated as follows: given three \emph{arbitrary balls} $\Omega_i=\B(a_i; r_i)$ for $i=1,2,3$ in the Euclidean plane, construct the smallest ball that covers all of the given balls. In this case, $I=\{1, 2, 3\}$, $J=\emptyset$, and problem (\ref{modeling}) reduces to
\begin{equation}\label{modeling1}
\mbox{\rm minimize }\mathcal{G}(x)=\max\{M(x; \Omega_1), M(x; \Omega_2), M(x; \Omega_3) \}, \; x\in \Bbb R^2.
\end{equation}
For two balls $\B(a; r)$ and $\B(b; s)$, we say that $\B(a; r)$ \emph{strictly contains} $\B(b; s)$ if $\B(b; s)\subseteq \B(a; r)$ and they have no boundary point in common, and $\B(a; r)$ \emph{tangentially contains} $\B(b; s)$ if $\B(b; s)\subseteq \B(a; r)$ and they have exactly one boundary point in common.

\begin{proposition}\label{p1} For $\ox\in \Bbb R^2$ and $r:=\mathcal{G}(\ox)$, one has that $|A(\ox)|=1$ and $\ox$ is the solution of problem {\rm(\ref{modeling1})} if and only if $\B(\ox; r)$ coincides with one of the ball $\B(a_i; r_i)$  for $i\in I$, and $\B(\ox; r)$ strictly contains the other balls.
\end{proposition}
\noindent {\bf Proof. }Observe that problem (\ref{modeling1}) has a unique solution by Theorem \ref{eu}. Suppose $|A(\ox)|=1$, say $A(\ox)=\{1\}$. Since $\ox$ is the unique solution of the optimization problem (\ref{modeling1}), one has
\begin{equation*}
0\in\partial M(\ox;\Omega_1).
\end{equation*}

Moreover, $r=\mathcal{G}(\ox)=M(\ox;\Omega_1)$, and $r>M(\ox; \Omega_i)$ for $i=2,3$. Since $0\in \partial M(\ox;\Omega_1)$, it follows from the subdifferential representation for $M(\cdot;\Omega_1)$ from Proposition \ref{subr} that $\ox=a_1$, and hence $$r=M(\ox;\Omega_1)=M(a_1; \Omega_1)=r_1.$$ Moreover, $r>M(a_1, \Omega_i)=\|a_1-a_i\|+r_i$ for $i=2,3$. This implies $\B(\ox; r)=\B(a_1, r_1)$ strictly contains two other balls.

The proof of the converse is straightforward. Indeed, assume that $\B(\ox; r)$ coincides with $\B(a_1; r_1)$ and strictly contains the other balls. Then $\ox=a_1$, $\mathcal{G}(\ox)=r=r_1=M(\ox; \Omega_1)$, and $\mathcal{G}(\ox)>\|a_1-a_i\|+r_i=M(\ox; \Omega_i)$ for $i=2,3$. Thus $A(\ox)=\{1\}$, and $0\in \partial \mathcal{G}(\ox)=\B.$ Therefore, $\ox$ is the solution of problem (\ref{modeling1}). $\h$ \vspace*{0.05in}

We are going to use the following notations when they are well-defined:
\begin{align*}
&u_1:=\mathcal P(a_2; \Omega_3), v_1:=\mathcal P(a_3; \Omega_2), w_1:=\dfrac{u_1+v_1}{2}, t_1:=\mathcal P(w_1; \Omega_1);\\
&u_2:=\mathcal P(a_3; \Omega_1), v_2:=\mathcal P(a_1; \Omega_3), w_2:=\dfrac{u_2+v_2}{2}, t_2:=\mathcal P(w_2; \Omega_2);\\
&u_3:=\mathcal P(a_1; \Omega_2), v_3:=\mathcal P(a_2; \Omega_1), w_3:=\dfrac{u_3+v_3}{2}, t_3:=\mathcal P(w_3; \Omega_3).
\end{align*}

\begin{figure}[!ht]
\begin{minipage}[b]{0.4\textwidth}
\centering
\includegraphics[width=5cm]{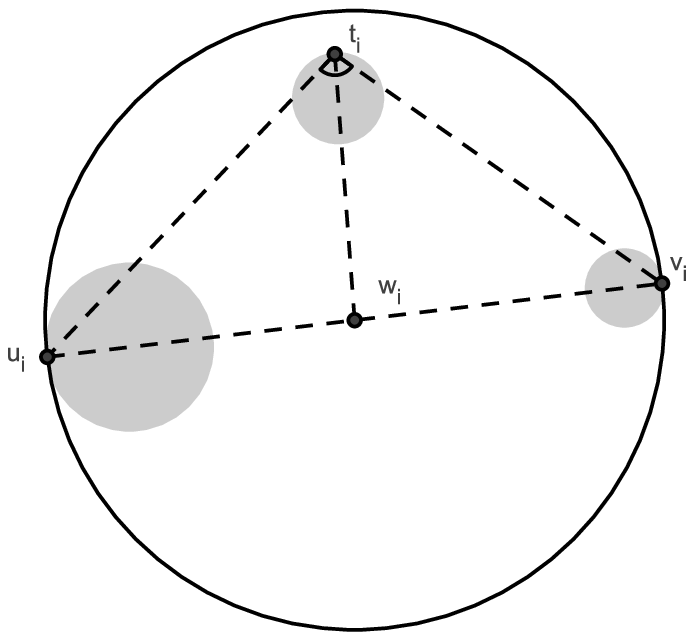}
\caption{$|A(\ox)|=2$}\label{fig:*1}
\end{minipage}
\hfill
\begin{minipage}[b]{0.4\textwidth}
\centering
\includegraphics[width=7cm]{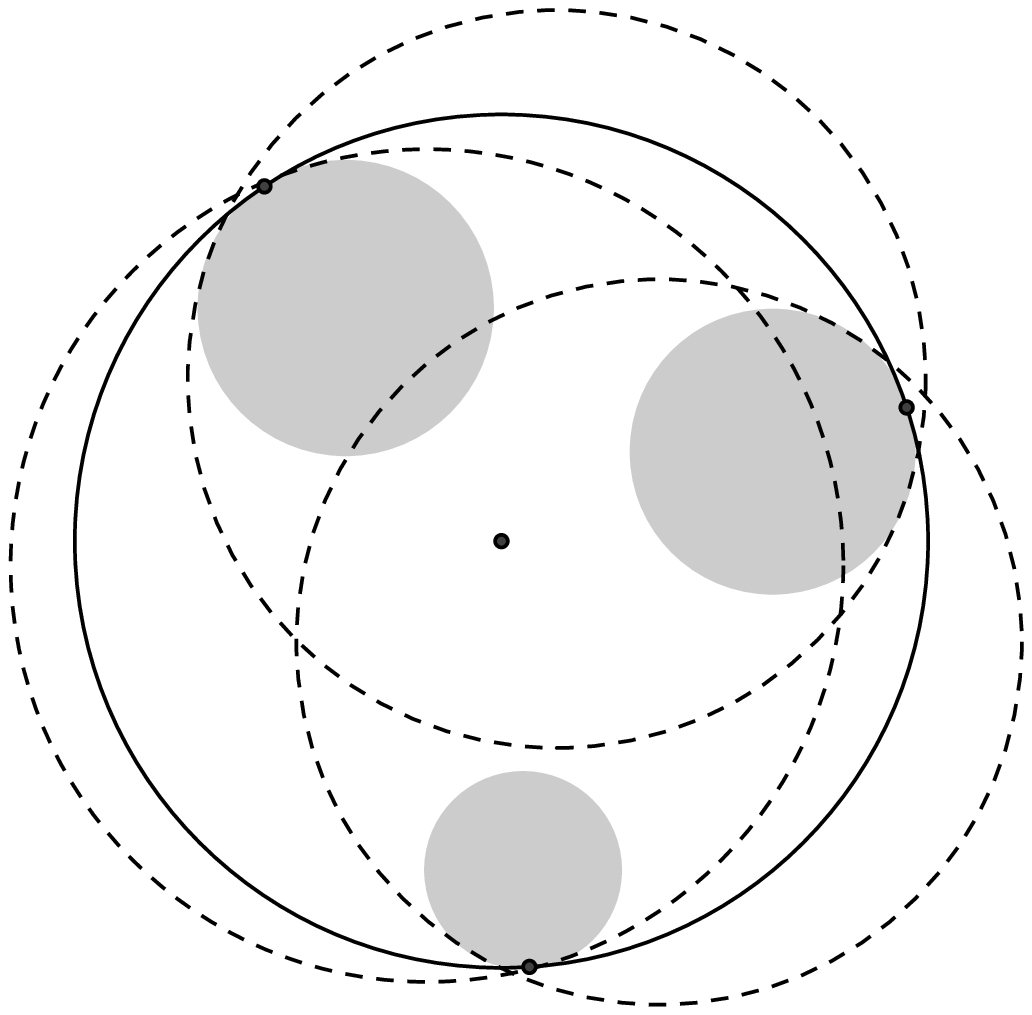}
\caption{$|A(\ox)|=3$}\label{fig:*2}
\end{minipage}
\end{figure}

For $a, b, c\in \Bbb R^2$, the notation $\widehat{bac}$ denotes the angle formed by vectors $\overrightarrow{ab}$ and $\overrightarrow{ac}$.

\begin{proposition} For $\ox\in \Bbb R^2$ and $r:=\mathcal{G}(\ox)$, one has that $|A(\ox)|=2$ and $\ox$ is the solution of problem {\rm(\ref{modeling1})} if and only if one of the following conditions holds:\\[1ex]
{\rm (1) } $\B(\ox; r)$ coincides with two balls among $\B(a_i; r_i)$  for $i\in I$ and strictly contains the remaining ball;\\
{\rm (2) } $\B(\ox; r)$ coincides with one of the balls among $\B(a_i; r_i)$  for $i\in I$, $\B(\ox; r)$ strictly contains another ball, and tangentially contains the remaining ball;\\
{\rm (3)} There exists $i\in I$ such that $\widehat{u_{i}t_{i}v_{i}}>90^{\circ}$, $\ox=w_i$, and $r=\dfrac{\|u_i-v_i\|}{2}$.
\end{proposition}
\noindent {\bf Proof.} Suppose that $|A(\ox)|=2$, for instance $A(\ox) = \{1, 2\}$. Consider the following cases:\\[1ex]
{\bf Case 1: }If both $\mathcal P(\ox; \Omega_1)$ and $\mathcal P(\ox; \Omega_2)$ are not singletons, then $\ox=a_1=a_2$. Moreover, $$r=M(\ox; \Omega_1)=M(\ox; \Omega_2)=r_1=r_2>M(\ox;\Omega_3),$$ so $\B(\ox; r)=\Omega_1=\Omega_2$, and $\B(\ox; r)$ strictly contains $\Omega_3$. In this case, (1) holds.\\[1ex]
{\bf Case 2: }If only one of the sets $\mathcal P(\ox; \Omega_1), \mathcal P(\ox; \Omega_2)$, for instance $\mathcal P(\ox; \Omega_1)$, is not a singleton, then one  has
$$\ox= a_1,  \|a_1-a_2\| = r_1 -r_2  \mbox{ and } \|a_1-a_3\| < r_1-r_3.$$
Thus, $\B(\ox; r)=\B(a_1; r_1)$, $\B(\ox; r)$ tangentially contains $\B(a_2; r_2)$, and $\B(\ox; r)$ strictly contains $\B(a_3; r_3)$. In this case, (2) holds.\\[1ex]
 {\bf Case 3: }If both $p_1:=\mathcal P(\ox; \Omega_1)$ and $p_2: = \mathcal P(\ox; \Omega_2)$ are singletons, one has by Theorem \ref{optimality} that $\ox=\dfrac{p_1+p_2}{2}$. Moreover,  $\|\ox-a_3\| +r_3< \dfrac{\|p_1-p_2\|}{2}$. In this case, using the representation for farthest projections, we have that $\ox$ belongs to the open line segment $(a_1, a_2)$ that connects $a_1$ and $a_2$, and  $\ox=w_3$. We also have that $p_1=u_3$, $p_2=v_3$, and $\widehat{u_3t_3v_3}$ is greater than $90^\circ$. In this case,  (3) holds.

The converse under (1) or (2) follows directly from Theorem \ref{optimality}. Suppose (3) holds. Then $$M(\ox;\Omega_1)=M(\ox;\Omega_2)>M(\ox;\Omega_3).$$ So $A(\ox)=\{1,2\}$. Since $\mathcal P(\ox;\Omega_1)=u_3$, $\mathcal P(\ox;\Omega_2)=v_3$, and $\ox=\dfrac{u_3+v_3}{2}$,  by Theorem \ref{optimality}, the element $\ox$ is the solution of problem (\ref{modeling1}).  $\h$

\begin{proposition}\label{ea} For $\ox\in \Bbb R^2$ and $r:=\mathcal{G}(\ox)$, one has that $|A(\ox)|=3$ and $\ox$ is the solution of problem {\rm(\ref{modeling1})} if and only if one of the following conditions holds:\\[1ex]
 {\rm (1)} $\B(\ox; r)$ coincides with one of the balls among $\Omega_i$ for $i\in I$ and tangentially contains two other balls.\\
 {\rm (2) } $\B(\ox; r)$ coincides with two of the balls among $\Omega_i$ for $i\in I$ and tangentially contains the remaining one.\\
 {\rm (3) } $\B(\ox; r)$ coincides with all three balls.\\
 {\rm (4)} $\mathcal P(\ox; \Omega_i)$ are singletons for $i\in I$, $\ox\in \mbox{\rm co }\{p_1, p_2, p_3\}$, where $p_i:=\mathcal P(\ox; \Omega_i)$, and $r=\|\ox-p_1\|=\|\ox-p_2\|=\|\ox-p_3\|.$
\end{proposition}
\noindent {\bf Proof. }Let us prove the implication {\it ``if''}. We have that
\begin{equation*}
0\in \partial \mathcal{G}(\ox)=\mbox{\rm co }\{\partial M(\ox; \Omega_1), \partial M(\ox; \Omega_2), \partial M(\ox; \Omega_3)\},
\end{equation*}
and $r=M(\ox; \Omega_1)=M(\ox; \Omega_2)=M(\ox; \Omega_3)$. Consider the following cases:\\[1ex]
{\bf Case 1: }At least one of the sets among $\mathcal P(\ox; \Omega_i)$ for $i\in I$ is not a singleton. If exactly one of the sets among $\mathcal P(\ox; \Omega_i)$ for $i\in I$ is not a singleton, say $\mathcal P(\ox;\Omega_1),$ then $\ox=a_1$ and $r=r_1$. Moreover, $$M(a_1; \Omega_i)=r=\|a_1-a_i\|+r_i \mbox{ for }i=2,3.$$ In this case, $\Omega_1$ tangentially contains two other balls. In the case exactly two sets among $\mathcal P(\ox; \Omega_i)$ for $i\in I$ are not singletons, then (2) holds, and if all $\mathcal P(\ox; \Omega_i)$ for $i\in I$ are not singletons, then (3) holds.\\[1ex]
{\bf Case 2: }All $\mathcal P(\ox; \Omega_i)$ for $i\in I$ are singletons. In this case, one has $\|\ox-p_i\|=r$ for $i\in I$ and by Theorem \ref{optimality}, $\ox\in \mbox{\rm co }\{p_1, p_2, p_3\}$.

The converse follows directly from Theorem \ref{optimality}. $\h$ \vspace*{0.05in}

We can construct the smallest enclosing ball of three given balls in the plane as follows:\\[1ex]
{\bf Step 1.} If there exists one of the three given balls that contains two remaining ones, for instance, $\Omega_2 \cup \Omega_3 \subseteq \Omega_1$, then $\ox=a_1$ is the solution of the problem and $r=r_1$ is the optimal value. Otherwise, go to next step.\\[1ex]
\noindent {\bf Step 2. }If one of three angles $\widehat{u_{i}t_{i}v_{i}}$, $i\in I$, is greater than $90^{\circ}$, then $w_i$ is the optimal solution and $r=\dfrac{\|u_i-v_i\|}{2}$ is the optimal value. Otherwise, go to next step.\\[1ex]
\noindent {\bf Step 3.} The smallest enclosing ball coincides with the Apollonius ball that is internally tangent to $\Omega_i$ for $i\in I$; see, e.g., \cite{gr} and the references therein.\vspace*{0.05in}

\noindent {\bf Three-Ball Problem: Model II.} The second model we consider in this subsection is: given three balls in $\Bbb R^2$ which are $\Omega_i=\B(a_i; r_i)$ for $i=1,2$ and $\Theta_1=\B(b_1; s_1)$, find the smallest ball that covers $\Omega_1$ and $\Omega_2$ and intersects $\Theta_1$. In this case, $I=\{1,2\}$, $J=\{1\}$, and problem (\ref{modeling}) reduces to
\begin{equation}\label{modeling3}
\mbox{\rm minimize }\mathcal{G}(x)=\max\{M(x; \Omega_1), M(x; \Omega_2), D(x; \Theta_1) \}, \; x\in \Bbb R^2.
\end{equation}
For any $u\in \Bbb R^2$, one has $|K(u)|\in \{0, 1, 2\}$, $|L(u)|\in \{0, 1\}$, and $1\leq |K(u)|+|L(u)|\leq 3$. In this case, problem (\ref{modeling3}) has a unique optimal solution by Theorem \ref{eu}.

We say that two balls \emph{strictly intersect} if they intersect at more than one points, and \emph{tangentially intersect} if they intersect each other at exactly one point.

\begin{proposition} For $\ox\in \Bbb R^2$ and $r:=\mathcal{G}(\ox)$, one has that $\ox$ is the solution of the problem {\rm (\ref{modeling3})} and \\$|K(\ox)|+|L(\ox)|=1$ if and only if $\B(\ox; r)$ coincides with one of the sets  $\Omega_i$ for i=1,2, strictly contains the other, and strictly intersects $\Theta_1$.
\end{proposition}
\noindent {\bf Proof. }Suppose that $|K(\ox)|+|L(\ox)|=1$. Then $|K(\ox)|=1$ and $|L(\ox)|=0$, or $|K(\ox)|=0$ and $|L(\ox)|=1$.\\[1ex]
{\bf Case 1: }$|K(\ox)|=1$ and $|L(\ox)|=0$. Suppose that $K(\ox)=\{1\}$. Then
\begin{equation*}
0\in \partial \mathcal{G}(\ox)=\partial M(\ox; \Omega_1).
\end{equation*}
This implies $\ox=a_1$, $\mathcal{G}(\ox)=M(\ox; \Omega_1)=r_1>M(\ox; r_2)=\|a_1-a_2\|+r_2$, and $$\mathcal{G}(\ox)=M(\ox; \Omega_1)=r_1>D(\ox; \Theta_1)=\|a_1-b_1\|-s_1.$$ In this case, we get the conclusion.\\[1ex]
{\bf Case 2: }$|K(\ox)|=0$ and $|L(\ox)|=1$. In this case, one has $\mathcal{G}(\ox)=D(\ox; \Theta_1)>0$ and
\begin{equation*}
0\in \partial \mathcal{G}(\ox)=\partial D(\ox; \Theta_1).
\end{equation*}
Then $\ox\notin \Theta_1$, and hence $\partial D(\ox; \Theta_1)=\dfrac{\ox-b_1}{\|\ox-b_1\|}$. We arrive at a contradiction since $\ox\neq b_1$. The converse is straightforward. $\h$ \vspace*{0.05in}

For $\{i, j\}=\{1, 2\}$, we will use the following notations when they are defined:
\begin{align*}
&u_1:=\mathcal P (b_1; \Omega_1), v_1:=\Pi (a_1; \Theta_1), w_1=\dfrac{u_1+v_1}{2}, t_1:=\mathcal P (w_1; \Omega_2);\\
&u_2:=\mathcal P (b_1; \Omega_2), v_2:=\Pi (a_2; \Theta_1), w_2=\dfrac{u_2+v_2}{2}, t_2:=\mathcal P (w_2; \Omega_1).
\end{align*}
\begin{proposition} For $\ox\in \Bbb R^2$ and $r:=\mathcal{G}(\ox)$, one has that $\ox$ is the solution of the problem {\rm (\ref{modeling3})} and $|K(\ox)|=|L(\ox)|=1$ if and only if one of the following conditions hold:\\[1ex]
{\rm (1)} $\B(\ox; r)$ coincides with one of the sets  $\Omega_i$ for i=1,2, strictly contains the other, and tangentially intersects $\Theta_1$.\\
{\rm (2)} One of the angles $\widehat{u_it_iv_i}$ for $i=1,2$ is greater than $90^{\circ}$, $\ox=w_i$, and $r=\dfrac{\|u_i-v_i\|}{2}$.
\end{proposition}
\noindent {\bf Proof. }Suppose that $\ox$ be a solution of the problem and $|K(\ox)|=|L(\ox)|=1$. Then $K(\ox)=\{1\}$ or $K(\ox)=\{2\}$ and $L(\ox)=\{1\}$. Consider the case where $K(\ox)=\{1\}$ and $L(\ox)=\{1\}$. Then $$r=\mathcal{G}(\ox)=M(\ox; \Omega_1)=D(\ox; \Theta_1)>M(\ox; \Omega_2)$$ and
\begin{equation*}
0\in \partial\mathcal{G}(\ox)=\mbox{\rm co }\{\partial M(\ox;\Omega_1), \partial D(\ox; \Theta_1)\}.
\end{equation*}
Since $D(\ox; \Theta_1)>0$, one has $\ox\notin\Theta_1$. Observe that if $\mathcal P (\ox;\Omega_1)$ is not a singleton, then $\ox=a_1$, and hence $r_1>M(\ox; \Omega_2)=\|a_1-a_2\|+r_2$, and $r_1=M(\ox;\Omega_1)=D(\ox; \Theta_1)=\|a_1-b_1\|-s_1$. In this case, (1) holds.

Suppose that $\mathcal P(\ox;\Omega_1)$ be a singleton. By Theorem \ref{optimality}, $\ox=\dfrac{y_1+z_1}{2}$, where $y_1=\mathcal P(\ox; \Omega_1)$ and \\$z_1=\Pi (\ox; \Theta_1)$. In this case, $y_1=u_1$, $z_1=v_1$, and $\ox=w_1$. Since $r=\dfrac{\|u_1-v_1\|}{2}>M(w_1;\Omega_2)=\|w_1-t_1\|$, one has that $\widehat{u_1t_1v_1}>90^\circ$. The proof of the converse is also straightforward. $\h$ \vspace*{0.05in}

Denote
\begin{align*}
x_1:=\mathcal P(a_2; \Omega_1), x_2:=\mathcal P(a_1;\Omega_2), y=\dfrac{x_1+x_2}{2}, z:=\Pi(y; \Theta_1).
\end{align*}
\begin{proposition} For $\ox\in \Bbb R^2$ and $r:=\mathcal{G}(\ox)$, one has that $\ox$ is the solution of the problem {\rm (\ref{modeling3})} and $|K(\ox)|=2$, $|L(\ox)|=0$ if and only if one of the following conditions hold:\\[1ex]
{\rm (1)} $\B(\ox; r)$ coincides with both $\Omega_i$ for $i=1,2$, and strictly intersects $\Theta_1$.\\
{\rm (2)} $\B(\ox; r)$ coincides with one of the balls $\Omega_i$ for $i=1,2$, tangentially contains the other, and strictly intersects $\Theta_1$.\\
{\rm (3)} The angle $\widehat{x_1zx_2}$ is greater than $90^\circ$, $\ox=y$, and $r=\dfrac{\|x_1-x_2\|}{2}.$
\end{proposition}
\noindent {\bf Proof. }Suppose that $\ox$ be a solution of problem (\ref{modeling3}), $K(\ox)=\{1,2\}$, and $L(\ox)=\emptyset$. Then $$r=\mathcal G(\ox) =M(\ox, \Omega_1)=M(\ox,\Omega_2) >D(\ox, \Theta_1),$$ and
\begin{equation*}
0\in \partial\mathcal{G}(\ox)=\mbox{\rm co }\{\partial M(\ox;\Omega_1), \partial M(\ox; \Omega_2)\}.
\end{equation*}

\begin{figure}[!ht]
\begin{minipage}[b]{0.4\textwidth}
\centering
\includegraphics[width=7.5cm]{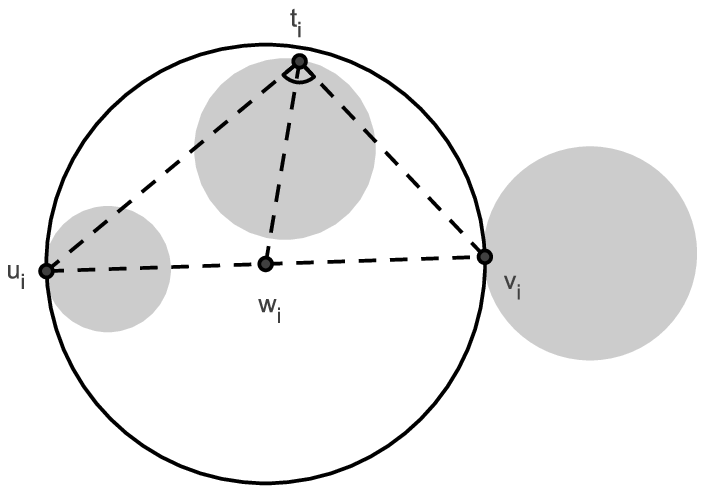}
\caption{$|K(\ox)|=1$, $|L(\ox)|=1$}\label{fig:*1}
\end{minipage}
\hfill
\begin{minipage}[b]{0.4\textwidth}
\centering
\includegraphics[width=6cm]{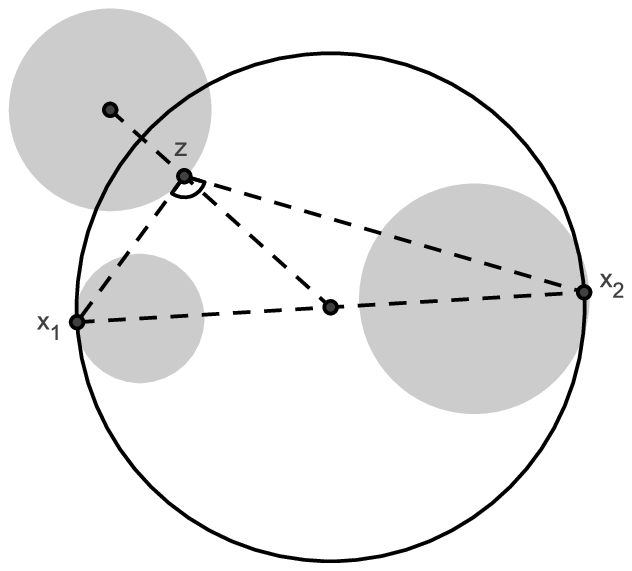}
\caption{$|K(\ox)|=2$, $|L(\ox)|=0$}\label{fig:*2}
\end{minipage}
\end{figure}

Let us consider the following cases:\\[1ex]
{\bf Case 1:} If $\ox=a_1$ and $\ox=a_2$, then $r_1=r_2=r, $ and hence $\B(\ox; r)=\Omega_1 =\Omega_2.$  On the other hand, since $D(\ox, \Theta_1)<r$ , one has that $\B(\ox; r)$ strictly intersects $\Theta_1.$\\
{\bf Case 2:} If $\ox=a_1$ and $\ox\ne a_2$, then $r_1=r =r_2+\|\ox -a_2\|.$ In this case, we have  $\B(\ox; r)=\Omega_1. $ Since $\|\ox-a_2\| =r_1-r_2$ and $D(\ox, \Theta_1)<r$, we have that $\B(\ox; r)$ tangentially contains $\Omega_2$ and strictly intersects $\Theta_1.$\\
{\bf Case 3:} Now, we consider the case  $\ox\ne a_1$ and $\ox\ne a_2.$  Using Theorem 3.2, we see that $$\ox =tx_1 +(1-t) x_2,\ t\in \ ]0,1[.$$ Since the problem has a unique solution, then $t=1/2$, i.e., $\ox =\dfrac{x_1+x_2}2=y$  and $r=\dfrac{\|x_1-x_2\|}2.$ Since $D(\ox, \Theta_1)<r,$ or $\|y-z\|<r,$ the angle
$\widehat{x_1zx_2}$ is greater than $90^\circ.$

The proof of the sufficient condition is straightforward. $\h$
\begin{proposition}\label{sa} For $\ox\in \Bbb R^2$ and $r:=\mathcal{G}(\ox)$, one has that $\ox$ is the solution of the problem {\rm (\ref{modeling3})} and $|K(\ox)|=2$, $|L(\ox)|=1$ if and only if one of the following conditions hold:\\[1ex]
{\rm (1)} $\B(\ox; r)$ coincides with both of the balls $\Omega_i$ for $i=1,2$, and tangentially intersects $\Theta_1$.\\
{\rm (2)} $\B(\ox; r)$ coincides with one of the balls $\Omega_i$ for $i=1,2$, tangentially contains the other, and tangentially intersects $\Theta_1$.\\
{\rm (3)} $\mathcal P(\ox; \Omega_i)$ for $i=1,2$ and $\Pi(\ox; \Theta_1)$ are singletons, $\|\ox-p_1\|=\|\ox-p_2\|=\|\ox-q_1\|$, where $p_i=\mathcal P(\ox; \Omega_i)$ for $i=1,2$ and $q_1:=\Pi(\ox; \Theta_1)$, and $\ox\in \mbox{\rm co }\{p_1, p_2, q_1\}$.
\end{proposition}
\noindent {\bf Proof. }Suppose that $\ox$ be a solution of the problem and $K(\ox)=\{1,2\}$ and $L(\ox)=\{1\}.$
 Then $$r=\mathcal G(\ox) =M(\ox, \Omega_1)=M(\ox,\Omega_2)=D(\ox, \Theta_1),$$ and
\begin{equation*}
0\in \partial\mathcal{G}(\ox)=\mbox{\rm co }\{\partial M(\ox;\Omega_1), \partial M(\ox; \Omega_2), \partial D(\ox, \Theta_1)\}.
\end{equation*}
We consider the following cases:\\[1ex]
{\bf Case 1:} If $\ox=a_1$ and $\ox=a_2$, then $r_1=r_2=r$, and hence $\B(\ox; r)=\Omega_1 =\Omega_2.$  On the other hand, since $D(\ox, \Theta_1)=r$, one has that $\B(\ox; r)$ tangentially intersects $\Theta_1.$\\
{\bf Case 2:} If $\ox=a_1$ and $\ox\ne a_2$, then $r_1=r =r_2+\|\ox -a_2\|.$ In this case, we have  $\B(\ox; r)=\Omega_1. $ Since $\|\ox-a_2\| =r_1-r_2$ and $D(\ox, \Theta_1)=r$, one has that $\B(\ox; r)$ tangentially contains $\Omega_2$ and tangentially intersects $\Theta_1$. \\
{\bf Case 3:} Now, we consider the case where  $\ox\ne a_1$ and $\ox\ne a_2.$  Then $\mathcal P(\ox; \Omega_i)$ for $i=1,2$ are singletons. Since $D(\ox, \Theta_1) =r>0$,  $\Pi (\ox, \Theta_1)$ is also a singleton.

Using Theorem \ref{optimality}, we see that $\ox =t_1p_1 +t_2p_2+t_3q_1, \ t_i\in [0,1]$, where $t_1+t_2+t_3=1.$ Moreover, $\|\ox-p_1\|=\|\ox-p_2\|=\|\ox-q_1\|.$

The proof of the sufficient condition is straightforward. $\h$ \vspace*{0.05in}

\begin{figure}[!ht]
\centering
\includegraphics[width=10cm]{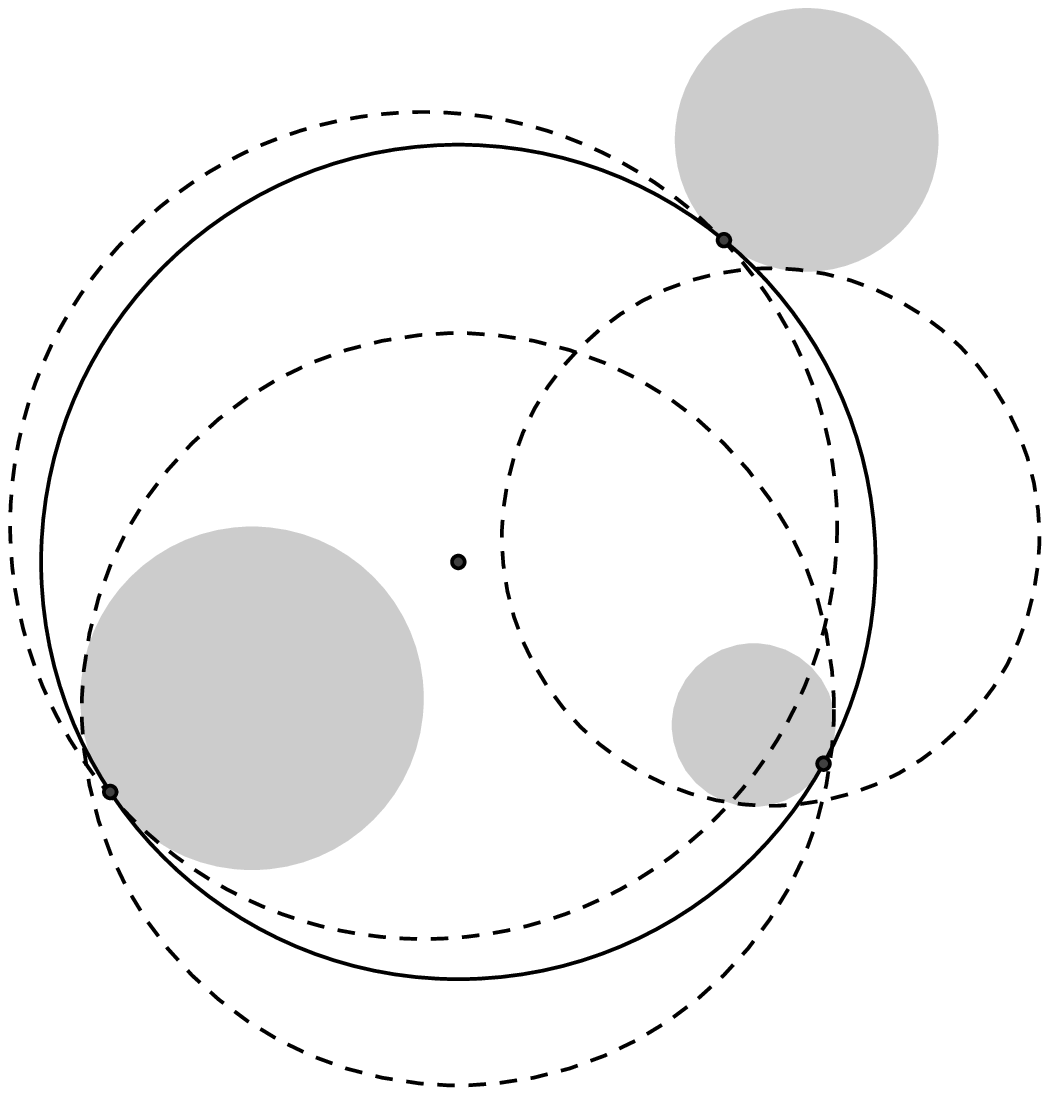}
\caption{$|K(\ox)|=2$, $|L(\ox)|=1$}\label{fig:}
\end{figure}

We are now able to construct the smallest ball that corresponds to the solution of problem (\ref{modeling3}) as follows:\\[1ex]
\noindent{\bf Step 1.} If there is a ball from $\Omega_1$ and $\Omega_2$ that contains the other and intersects $\Theta_1$, then that ball is the solution of the generalized Sylvester problem. Otherwise, we go to next step.\\
\noindent{\bf Step 2.} One of the angles $\widehat{u_it_iv_i}$ for $i=1,2$ is greater than $90^{\circ}$, $\ox=w_i$, and $r=\dfrac{\|u_i-v_i\|}{2}$. Otherwise, we go to next step.\\
\noindent {\bf Step 3.} The angle $\widehat{x_1zx_2}$ is greater than $90^\circ$, $\ox=y$, and $r=\dfrac{\|x_1-x_2\|}{2}.$ Otherwise, we go to next step.\\
\noindent{\bf Step 4.} In this case, the smallest ball is the Apollonius ball that is internally tangent to $\Omega_1, \Omega_2$ and externally tangent to $\Theta$. \vspace*{0.05in}

\noindent {\bf Three-Ball Problem: Model III.} The third model we will consider is: given three  balls $\Omega_1=\B(a_1; r_1)$,  $\Theta_1=\B(b_1; s_1)$, and $ \Theta_2=\B(b_2; s_2)$, find the smallest ball that covers $\Omega_1$ and intersects $\Theta_1$ and $\Theta_2$.  In this case, $I=\{1\}$, $J=\{1,2\}$, and problem (\ref{modeling}) reduces to:
\begin{equation}\label{modeling2}
\mbox{\rm minimize }\mathcal{G}(x)=\max\{M(x; \Omega_1), D(x; \Theta_1), D(x; \Theta_2) \}, \; x\in \Bbb R^2.
\end{equation}

\begin{proposition} For $\ox\in \Bbb R^2$ and $r:=\mathcal{G}(\ox)$, one has that $\ox$ is the solution of the problem {\rm (\ref{modeling2})} and $|K(\ox)|+|L(\ox)|=1$ if and only if $\B(\ox; r)$ coincides with $\Omega_1$ and strictly intersects both $\Theta_1$ and $\Theta_2$.
\end{proposition}
\noindent {\bf Proof. }Suppose that $|K(\ox)|+|L(\ox)|=1$. Then $|K(\ox)|=1$ and $|L(\ox)|=0$, or $|K(\ox)|=0$ and $|L(\ox)|=1$. Consider the first case where $|K(\ox)|=1$ and $|L(\ox)|=0$. Then $K(\ox)=\{1\}$ and
\begin{equation*}
0\in \partial M(\ox; \Omega_1).
\end{equation*}
Moreover, $r=M(\ox; \Omega_1)>D(\ox; \Theta_i)$ for $i=1,2$. This implies $\ox=a_1$, and hence $$M(\ox;\Omega_1)=r_1=r>D(\ox; \Theta_i).$$ In this case, $\B(\ox; r)$ coincides with $\Omega_1$ and strictly intersects $\Theta_1$ and $\Theta_2$.

For the case where $|K(\ox)|=0$ and $|L(\ox)|=1$, we can assume that $L(\ox)=\{1\}$. Then
\begin{equation*}
0\in \partial D(\ox; \Theta_1),
\end{equation*}
$D(\ox; \Theta_1)>D(\ox; \Theta_2)$, and $D(\ox; \Theta_1)>M(\ox; \Omega_1)$. This implies $D(\ox; \Theta_1)>0$, and hence $\ox\notin\Theta_1$. We have arrived at a contradiction due to the representation of the subdifferential of distance function at out-of-set points from Proposition \ref{subr}. The converse follows from Theorem \ref{optimality}. $\h$ \vspace*{0.05in}

For $\{i, j\}=\{1,2\}$, we will use the following notations when they are defined:
$$c_i=\mathcal P(b_i;\Omega_1), d_i=\Pi (a_1; \Theta_i), e_i=\dfrac{c_i+d_i}{2}, f_i=\Pi (e_i; \Theta_j).$$
\begin{proposition} For $\ox\in \Bbb R^2$ and $r:=\mathcal{G}(\ox)$, one has that $\ox$ is the solution of problem {\rm (\ref{modeling2})} and $|K(\ox)|=|L(\ox)|=1$ if and only if one the following holds:\\[1ex]
{\rm (1) } $\B(\ox; r)$ coincides with $\Omega_1$, tangentially intersects one of the $\Theta_i$ for $i=1,2$ and strictly intersects the other.\\
{\rm (2) } There exists $i\in J$ such that $\widehat{c_if_id_i}>90^{\circ}$, $\ox=e_i$, and $r=\dfrac{\|c_i-d_i\|}{2}$.
\end{proposition}
\noindent {\bf Proof. }Suppose that $\ox$ be a solution of the problem and $|K(\ox)|=|L(\ox)|=1$. Then $K(\ox)=\{1\}$ and [$L(\ox)=\{1\}$ or $L(\ox)=\{2\}$]. Consider the case where $K(\ox)=\{1\}$ and $L(\ox)=\{1\}$. Then $$r=\mathcal{G}(\ox)=M(\ox; \Omega_1)=D(\ox; \Theta_1)>D(\ox; \Theta_2)$$ and
\begin{equation*}
0\in \partial\mathcal{G}(\ox)=\mbox{\rm co }\{\partial M(\ox;\Omega_1), \partial D(\ox; \Theta_1)\}.
\end{equation*}
Since $D(\ox; \Theta_1)>0$, one has $\ox\notin\Theta_1$.

In the case $\mathcal P (\ox;\Omega_1)$ is not a singleton, one has $\ox=a_1$. Then $r_1=M(\ox; \Omega_1)=D(\ox; \Theta_1)=\|a_1-b_1\|-s_1$ and $r_1=M(\ox; \Omega_1)>D(\ox; \Theta_2)$. In this case, (1) holds.

Now assume that $\mathcal P(\ox;\Omega_1)$ is a singleton. By Theorem \ref{optimality}, one has $\ox=\dfrac{z_1+y_1}{2}$, where $y_1=\mathcal P (\ox; \Omega_1)$ and $z_1={\Pi}(\ox; \Theta_1)$. In this case, $z_1=c_1$, $y_1=d_1$, and $\ox=e_1$. Since $r=\|c_1-d_1\|/2>D(e_1;\Theta_2)=\|e_1-f_1\|$, one has that $\widehat{c_if_id_i}>90^\circ$ for $i=1$. The converse follows from Theorem \ref{optimality}. $\h$ \vspace*{0.05in}

\begin{figure}[!ht]
\begin{minipage}[b]{0.4\textwidth}
\centering
\includegraphics[width=6cm]{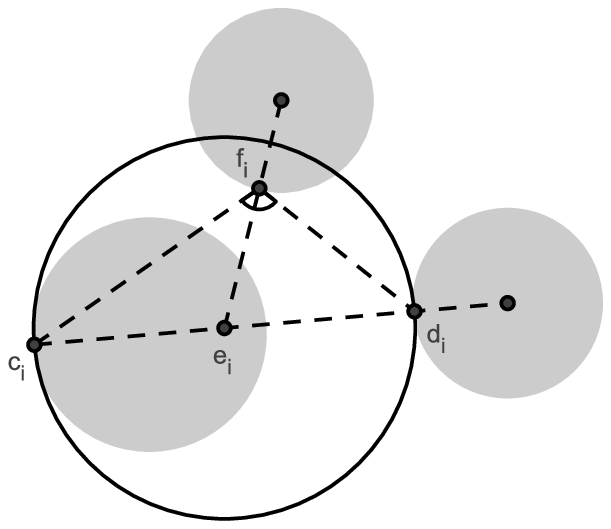}
\caption{$|K(\ox)|=1$, $|L(\ox)|=1$}\label{fig:*1}
\end{minipage}
\hfill
\begin{minipage}[b]{0.4\textwidth}
\centering
\includegraphics[width=8cm]{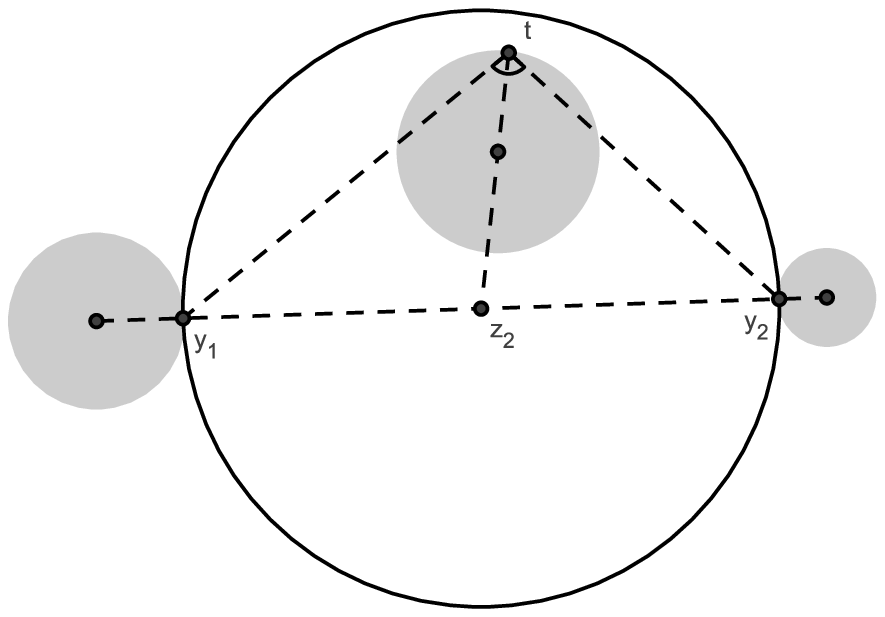}
\caption{$|K(\ox)|=0$, $|L(\ox)|=2$}\label{fig:*2}
\end{minipage}
\end{figure}

For $i, j\in \{1, 2\}$, $i\neq j$, we will use the following notations when they are defined:
$$y_j=\Pi (b_i, \Theta_j), \ z=\dfrac{y_1+y_2}{2},\  t=\mathcal P (z; \Omega_1).$$

\begin{proposition} For $\ox\in \Bbb R^2$ and $r:=\mathcal{G}(\ox)$, one has that $\ox$ is the solution of the problem {\rm (\ref{modeling2})} and $|K(\ox)|=0$, $|L(\ox)|=2$ if and only if  $\widehat{y_1t y_2}>90^{\circ}$, $\ox=z$, and $r=\dfrac{\|y_1-y_2\|}{2}$.
\end{proposition}
\noindent {\bf Proof. }In this case, $K(\ox)=\emptyset$ and $L(\ox)=\{1, 2\}$. Then
\begin{equation*}
0\in \partial \mathcal{G}(\ox)=\mbox{\rm co }\{\partial D(\ox; \Theta_1), \partial D(\ox; \Theta_2)\},
\end{equation*}
and $r=D(\ox; \Theta_1)=D(\ox; \Theta_2)>M(\ox; \Omega_1)$. This implies $\ox\notin \Theta_i$ for $i=1,2$. Let $z_i=\Pi (\ox; \Theta_i)$ for $i=1,2$. By Theorem \ref{optimality}, $\ox=\dfrac{z_1+z_2}{2}$. This implies $z_i=y_i$ for $i=1,2$ and $\ox=z$. Moreover, $$\dfrac{\|y_1-y_2\|} 2 > M(\ox; \Omega_1)=\|a_1-z\|.$$ Thus, $\widehat{y_1t y_2}>90^{\circ}$. The proof of the converse is straightforward. $\h$

\begin{figure}[!ht]
\centering
\includegraphics[width=10cm]{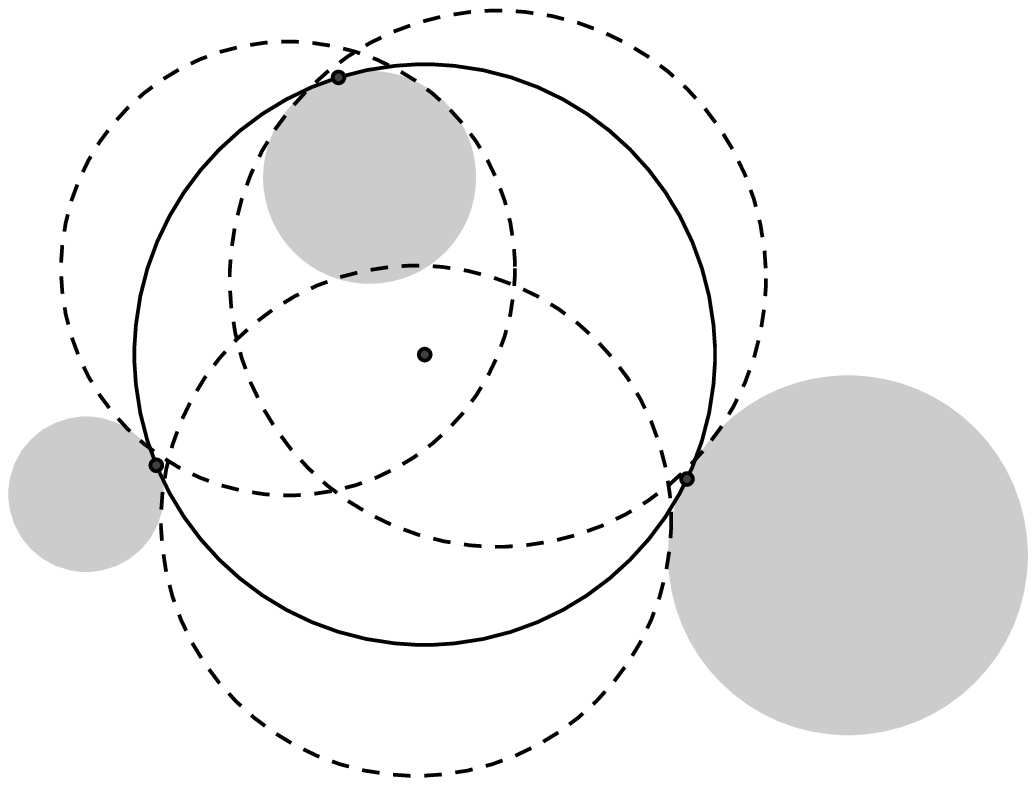}
\caption{$|K(\ox)|=1$, $|L(\ox)|=2$}\label{fig:*2}
\end{figure}

The proof of the proposition below is similar to that of Proposition \ref{sa}.
\begin{proposition} For $\ox\in \Bbb R^2$ and $r:=\mathcal{G}(\ox)$, one has that $\ox$ is the solution of the problem {\rm (\ref{modeling2})} and $|K(\ox)|=1$ and $L(\ox)|=2$ if and only if one the following holds:\\[1ex]
{\rm (1) } $\B(\ox; r)$ coincides with $\Omega_1$, is tangent (externally) to both $\Theta_i$ for $i=1,2$.\\
{\rm (2) } $\mathcal P (\ox; \Omega_1)$, $\Pi (\ox; \Theta_i)$ for $i=1,2$ are singletons, $\|\ox-p_1\|=\|\ox-q_1\|=\|\ox-q_2\|$, where $p_1=\mathcal P (\ox; \Omega_1)$ and $q_i:=\Pi (\ox; \Theta_i)$ for $i=1,2$, and $\ox\in \mbox{\rm co }\{p_1, q_1, q_2\}$.
\end{proposition}

We are now able to construct the smallest ball that corresponds to the solution of problem (\ref{modeling2}) as follows:

\noindent {\bf Step 1.} If $\Omega_1\cap \Theta_1 \neq \emptyset$ and $\Omega_1\cap \Theta_2 \neq \emptyset$ then $\ox=a_1$ is the solution of the problem and $r_1$ is the optimal value. Otherwise, go to next step.\\
\noindent {\bf Step 2.} If one of the angles $\widehat{c_if_id_i}$ for $i=1,2$ is greater than $90^{\circ}$, then $\ox=e_i$, $r=\dfrac{\|c_i-d_i\|}{2}$. Otherwise, go to the next step.\\
\noindent {\bf Step 3.} If the angle $\widehat{y_1t y_2}>90^{\circ}$, then $\ox=z$ and $r=\dfrac{\|y_1-y_2\|}{2}$. Otherwise, go to the next step.\\
\noindent{\bf Step 4.} In this case, the smallest ball is the Apollonius ball that is internally tangent to $\Omega_1$ and external tangent to $\Theta_1, \Theta_2$.\vspace*{0.05in}

\noindent {\bf Three-Ball Problem: Model IV.} Let us now consider the \emph{smallest intersecting ball problem}: given three balls $\Theta_i=\B(b_i; s_i)$ for $i=1,2,3$, find the smallest ball that intersects $\Theta_i$ for $i=1,2,3$. In this case, $I=\emptyset$, $J=\{1,2,3\}$, and problem (\ref{modeling}) reduces to
\begin{equation}\label{modeling4}
\mbox{\rm minimize }\mathcal{G}(x)=\max\{D(x; \Theta_1), D(x; \Theta_2), D(x; \Theta_3) \}, \; x\in \Bbb R^2.
\end{equation}
In the case $\cap_{i=1}^3\Theta_i\neq\emptyset$, any point in this intersection is a solution of problem (\ref{modeling4}), so we only consider the case where this intersection is empty. By Theorem \ref{eu}, problem (\ref{modeling4}) has a unique solution. It is also not hard to see that $|A(\ox)|=|L(\ox)|\geq 2$.\vspace*{0.05in}

We will use the following notations when they are defined:
\begin{align*}
&u_1=[b_2,b_3]\cap bd(\Theta_2),\
v_1=[b_2,b_3]\cap bd(\Theta_3),\notag\\
&u_2=[b_1, b_3]\cap bd(\Theta_3),\
v_2=[b_1, b_3]\cap bd(\Theta_1), \notag\\
&u_3=[b_1, b_2] \cap bd(\Theta_1),\
v_3=[b_1, b_2]\cap bd(\Theta_2),\notag\\
&m_1=\frac{u_1+v_1}{2},\ m_2=\frac{u_2+v_2}{2},\ m_3=\frac{u_3+v_3}{2},\notag \\
&x_1=\Pi(m_1;\Theta_1),\
x_2=\Pi(m_2; \Theta_2),\
x_3=\Pi(m_3;\Theta_3),\notag\\
&\B_1=\B(m_1; \frac{\chuan{u_1-v_1}}{2}), \ \B_2=\B(m_2; \frac{\chuan{u_2-v_2}}{2}),\  \B_3=\B(m_3; \frac{\chuan{u_3-v_3}}{2}). \label{datcacB}
\end{align*}

\begin{figure}[!ht]
\begin{minipage}[b]{0.4\textwidth}
\centering
\includegraphics[width=8cm]{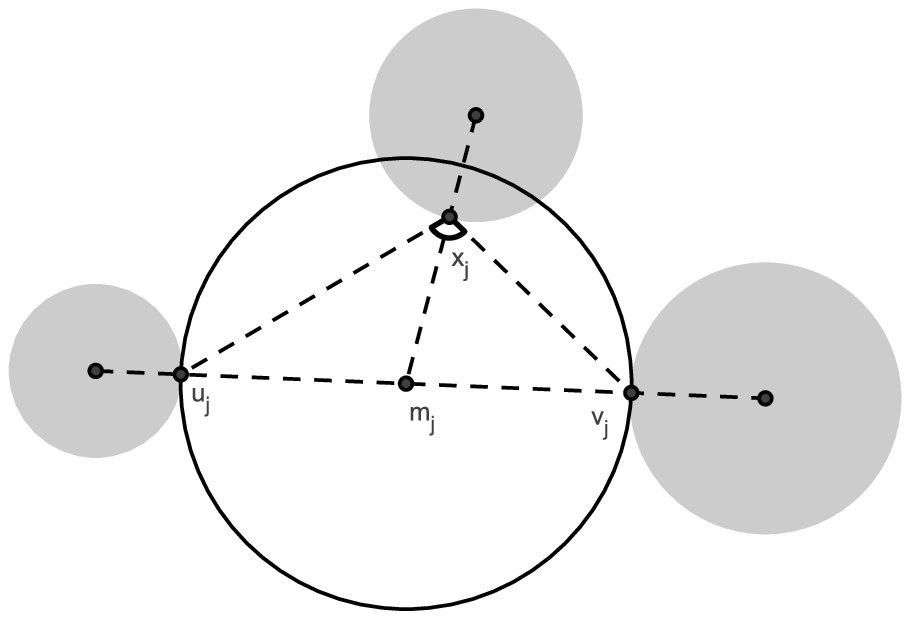}
\caption{$|A(\ox)|=2$}\label{fig:*1}
\end{minipage}
\hfill
\begin{minipage}[b]{0.4\textwidth}
\centering
\includegraphics[width=8cm]{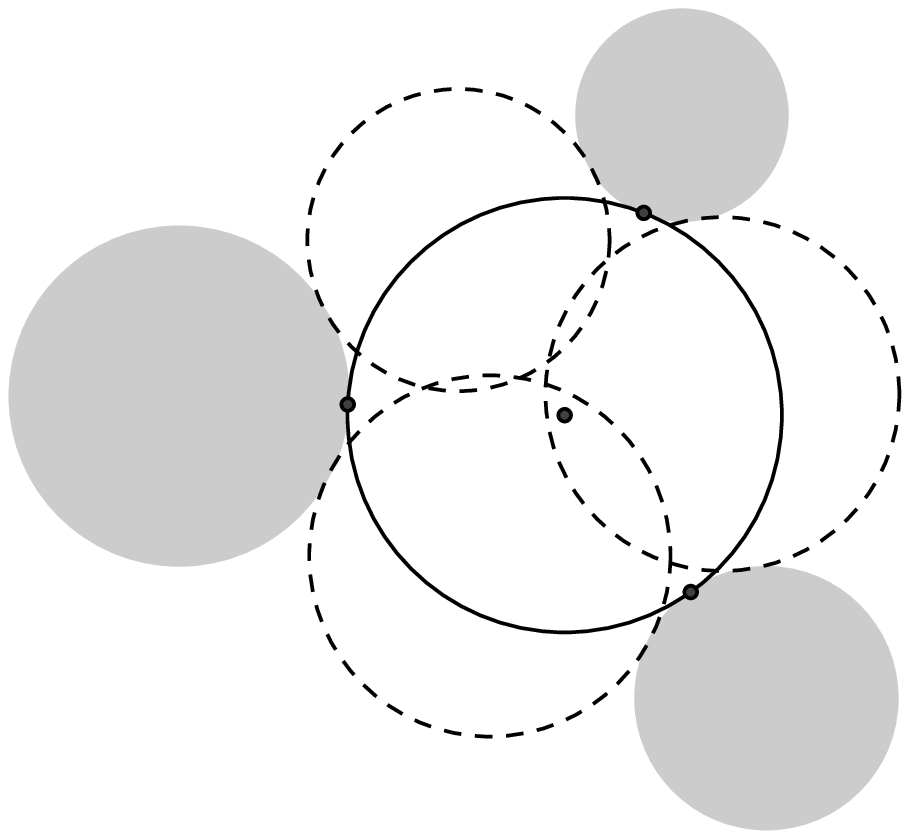}
\caption{$|A(\ox)|=3$}\label{fig:*2}
\end{minipage}
\end{figure}

It is not hard to see that the smallest intersecting ball can be constructed as below: \\[1ex]
{\bf Step 1:} If there exists $j\in \{1, 2, 3\}$ such that $\widehat{u_jx_jv_j}$ is greater than $90^\circ$, then $\B_j$ is the solution of the problem. Otherwise, go to the next step.\\[1ex]
{\bf Step 2:} The smallest intersecting ball is the Apollonius ball that is externally tangent to the three balls.

The readers are referred to \cite{naj} for a similar construction for the case where three given balls are disjoint.
\section{Generalized Fermat-Torricelli Problem for Euclidean Balls}
\label{sec:3}

In this section, we consider the following problem:
\begin{equation}\label{ft}
\mbox{\rm minimize }\mathcal{H}(x)=\sum_{i\in I}M(x; \Omega_i)+\sum_{j\in J}D(x;\Theta_j),\  x\in \Bbb R^2,
\end{equation}
where $|I|+|J|=3$. It is clear that $\ox$ is an optimal solution of problem (\ref{ft}) if and only if it is a solution of the following problem:
\begin{equation*}\label{ft1}
\mbox{\rm minimize }\mathcal{H}(x)=\sum_{i\in I}\|x-a_i\|+\sum_{j\in J}D(x;\Theta_j),\  x\in \Bbb R^2.
\end{equation*}
Since a singleton is a ball with radius $0$, it suffices to consider the following problem:
\begin{equation}\label{ft2}
\mbox{\rm minimize }\mathcal{H}(x)=\sum_{i=1}^3D(x;\Theta_i),\  x\in \Bbb R^2.
\end{equation}
Throughout this section, we assume that the given balls have \emph{distinct centers} since the other case is trivial.

\subsection{Existence and Uniqueness of Optimal Solutions}

In this subsection, we will study properties of solutions of problem (\ref{ft2}) that enable us to derive a necessary and sufficient condition for the problem to have a unique solution.
\begin{proposition} \label{pro1} The solution set $S$ of problem {\rm (\ref{ft2})} is a nonempty compact convex set in $\Bbb R^2$.
\end{proposition}
\noindent {\bf Proof. }It is clear that the function $\mathcal{H}(x)$ is continuous. For any $\gamma\geq\inf_{x\in \Bbb R^2}\mathcal{H}(x)$, we see easily that $\{x\in \Bbb R^2\; |\; \mathcal{H}(x)\leq\gamma\}$ is compact. Thus, {\rm (\ref{ft2}) always has an optimal solution. In particular, $S$ is compact. Since $\mathcal{H}$ is convex and continuous, the solution set of the problem is also convex. $\h$

For any $u\in \Bbb R^2$, define $$A(u):=\{i\in \{1, 2, 3\}\; |\; u\in \Theta_i\}.$$

\begin{proposition} Suppose that $|A(x_*)|=0$. Then $x_*$ is an optimal solution of problem {\rm (\ref{ft2})} if and only if $x_*$ is the solution of the classical Fermat-Torricelli generated by the centers of the balls: $b_1, b_2, b_3$.
\end{proposition}
\noindent {\bf Proof. }Suppose $|A(x_*)|=0$ and $x_*$ be an optimal solution of problem {\rm (\ref{ft2})}. Then $x_*\notin\Theta_i$ for $i=1,2,3$. Choose $\delta>0$ such that $\B(x_*; \delta)\cap \Theta_i=\emptyset$ for $i=1,2,3$. Then the following holds for every $x\in \B(x_*; \delta)$:
\begin{equation*}
\sum_{i=1}^3\|x_*-b_i\|-\sum_{i=1}^3s_i=\inf_{x\in \R^2}\mathcal{H}(x)=\mathcal{H}(x_*)\leq \mathcal{H}(x)=\sum_{i=1}^3\|x-b_i\|-\sum_{i=1}^3s_i.
\end{equation*}
Thus, $x_*$ is a local minimum of the problem
\begin{equation*}
\mbox{\rm minimize }\mathcal{F}(x):=\sum_{i=1}^3\|x-b_i\|, x\in \R^2,
\end{equation*}
so it is also an absolute minimum of the problem since $\mathcal{F}$ is a convex function. Therefore, $x_*$ the solution of the classical Fermat-Torricelli generated by $b_1, b_2, b_3$. The converse is also straightforward. $\h$

\begin{lemma}\label{lm1} Let $S$ be the solution set of problem {\rm (\ref{ft2})}. Suppose that there exists $x_*\in S$ such that $A(x_*)=\emptyset$. Then $S$ is a singleton, namely $S=\{x_*\}$.
\end{lemma}
\noindent {\bf Proof.} Suppose that $A(x_*)=\emptyset$, and there exists $y_*\in S$ with $x_*\neq y_*$. Since $A(x_*)=\emptyset$, one has $x_*\notin\Theta_i$ for $i=1,2,3$. Then $x_*$ is the solution of the classical Fermat-Torricelli problem generated by the centers of the balls: $b_1, b_2,$ and $b_3$. Since $S$ is convex, one has $[x_*, y_*]\subseteq S$. Thus, it is possible to find $z_*\neq x_*$, $z_*\in S$, and $z_*\notin\Theta_i$ for $i=1,2,3$. Then $z_*$ is also the solution of the classical Fermat-Torricelli problem generated by the centers of the balls, which is a contradiction since this problem has a unique solution. $\h$

\begin{lemma}\label{st}For any $\alpha>0$, and $a, b\in \Bbb R^2$ with $a\neq b$, consider the set
\begin{equation*}\label{Elip}
E:=\{x\in \Bbb R^2\; \; |\; \; \|x-a\| +\|x-b\|=\alpha\}.
\end{equation*}
Suppose that $\|a-b\|<\alpha$. For $x, y\in E$ with $x\neq y$, one has $\|\dfrac{x+y}{2}-a\|+\|\dfrac{x+y}{2}-b\|<\alpha$, and in particular, $\dfrac{x+y}{2}\notin E$.
\end{lemma}
\noindent {\bf Proof.} We have
$$\aligned
\|\dfrac{x+y}2 -a\| +\|\dfrac{x+y}2-b\| =& \dfrac 12(\|x-a+y-a\| +\|x-b +y-b\|)\\
\le &\dfrac 12 (\|x-a\| +\|y-a\| +\|x-b\| +\|y-b\| )=\alpha.
\endaligned$$
Suppose by contradiction that
\begin{equation}\label{contradiction}
\|\dfrac{x+y}2 -a\| +\|\dfrac{x+y}2-b\| =\alpha.
\end{equation}
That means $\dfrac{x+y}{2}\in E$. By a property of the Euclidean norm, one has
$$x-a= k(y-a)\quad \text {and}\quad x-b =m(y-b)$$
for some numbers $k, m$ in $(0, +\infty)\setminus \{1\}$ since $x\ne y $. This implies
\begin{equation*}
a=\dfrac{1}{1-k}x-\dfrac{k}{1-k}y\in L(x,y) \mbox{ and }b=\dfrac{1}{1-m}x-\dfrac{m}{1-m}y\in L(x,y).
\end{equation*}
Since $\alpha >\|a-b\|$, it is not hard to see that $x,y\in L(a,b)\setminus [a,b]$. One also has $\dfrac{x+y}2\in [a,b]$. Indeed, assume for instance, that the order of points is: $x$, $\dfrac{x+y}{2}$, $a, b$, $y$. Then
\begin{equation*}
\|\dfrac{x+y}{2}-a\| +\| \dfrac{x+y}{2}-b\| < \|x-a\|+\|x-b\|=\alpha,
\end{equation*}
which contradicts (\ref{contradiction}). Now, since $\dfrac{x+y}2\in [a,b]$, one has
\begin{equation*}
\|\dfrac{x+y}2 -a\| +\|\dfrac{x+y}2-b\|=\|a-b\|<\alpha,
\end{equation*}
which is a contradiction. We have proved that $\|\dfrac{x+y}2 -a\| +\|\dfrac{x+y}2-b\| <\alpha.$ Thus, $\dfrac {x+y}2\notin E.$ $\h$ \vspace*{0.05in}

Let $\Omega$ be a subset of $\R^n$. We say that $\Omega$ is \emph{strictly convex} if for any $x, y\in \Omega$ with $x\neq y$ and for any $t\in \ ]0,1[$, one has $tx+(1-t)y\in \mbox{\rm int }\Omega$. In the setting of Lemma \ref{st}, the set $$\Tilde{E}:=\{x\in \Bbb R^2\; \; |\; \; \|x-a\| +\|x-b\|\leq\alpha\}$$ is strictly convex. This is in fact a well-known result, but we provide the detailed proof for the convenience of the readers.

\begin{lemma}\label{lm2} Let $S$ be the solution set of problem {\rm (\ref{ft2})}. Suppose that there exists $x_*\in S$ such that $A(x_*)=\{i\}$ and $x_*\notin [b_j, b_k]$, where $i, j, k$ are distinct indices in $\{1, 2, 3\}$. Then $S$ is a singleton, namely $S=\{x_*\}$.
\end{lemma}
\noindent {\bf Proof. }Suppose without loss of generality that $A(x_*)=\{1\}$. Then and $x_*\in\Theta_1$, $x_*\notin\Theta_2$, $x_*\notin\Theta_3$, and $x_*\notin [b_2, b_3]$. Let $$\alpha:=\|x_*-b_2\|+\|x_*-b_3\|=\inf_{x\in \Bbb R^2}\mathcal{H}(x)+s_2+s_3.$$ Then $\alpha > \|b_2-b_3\|.$ Define
\begin{equation*}
E:=\{x\in \Bbb R^2\; |\; \|x-b_2\|+\|x-b_3\|=\alpha\}.
\end{equation*}
Clearly, $x_*\in E\cap\Theta_1$. Suppose by contradiction that $S$ be not a singleton. Then there exists $y_*\in S$ and $y_*\neq x_*$. Since $S$ is convex, $[x_*, y_*]\subseteq S$.  We can choose $z_*\in [x_*, y_*]$ that is close enough and distinct from $x_*$ such that $[x_*, z_*]\cap \Theta_2=\emptyset$ and $[x_*, z_*]\cap\Theta_3=\emptyset$.
Let us first show that $z_*\notin\Theta_1$. Indeed, assume that $z_*\in \Theta_1$. Then $D(z_*;\Theta_1)=0$, and hence $$\mathcal{H}(z_*)=\inf_{x\in \R^2} \mathcal{H}(x)=D(z_*;\Theta_2)+D(z_*; \Theta_3)=\|z_*-b_2\|-s_2+\|z_*-b_3\|-s_3.$$
This implies $z_*\in E$. Since $S$ and $\Theta_1$ is convex, $\dfrac{x_*+z_*}{2}\in S\cap \Theta_1$. Moreover, $\dfrac{x_*+z_*}{2}\notin \Theta_2$ and $\dfrac{x_*+z_*}{2}\notin\Theta_3$. Again, we have that $\dfrac{x_*+z_*}{2}\in E$. This is not the case by Lemma \ref{st}. We have shown that $z_*\in S$ and $A(z_*)=\emptyset$, that contradicts the result from Lemma \ref{lm1}. Therefore, $S$ must be a singleton.$\h$

\begin{lemma}\label{lm3} Let $S$ be the solution set of problem {\rm (\ref{ft2})}. Suppose that $x_*\in S$ and $A(x_*)=\{i, j\}\subset \{1,2,3\}$. Then $x_*\in \mbox{\rm bd }(\Theta_i\cap \Theta_j)$.
\end{lemma}
\noindent {\bf Proof. }Assume without loss of generality that $A(x_*)=\{1,2\}$. Then $x_*\in \Theta_1\cap \Theta_2$ and $x_*\notin\Theta_3$. Suppose by contradiction that $x_*\in \mbox{\rm int }(\Theta_1\cap \Theta_2)$. Then there exists $\delta>0$ with $\B(x_*; \delta)\subseteq \Theta_1\cap \Theta_2$. Let $p:=\Pi(x_*; \Theta_3)$ and $\gamma:=\|x_*-p\|>0$. Let $q\in [x_*, p]\cap \B(x_*, \delta)$ satisfy $\|q-p\|<\|x_*-p\|=D(x_*; \Theta_3)$. Then
\begin{equation*}
\mathcal{H}(q)=\sum_{i=1}^3 D(q; \Theta_i)=D(q; \Theta_3)\leq \|q-p\|<\|x_*-p\|=D(x_*; \Theta_3)=\sum_{i=1}^3D(x_*; \Theta_i)=\mathcal{H}(x_*).
\end{equation*}
This is a contradiction to the fact that $x_*\in S$. $\h$

\begin{theorem}\label{tmu} Suppose $\cap_{i=1}^3\Theta_i=\emptyset$. Problem {\rm (\ref{ft2})} has more than one solution if and only if there is a set $[b_i, b_j]\cap \Theta_k$, for distinct indices $i, j, k\in \{1, 2, 3\}$, that contains a point $u$ belonging to the interior of $\Theta_k$ and $|A(u)|=1$ .
\end{theorem}
\noindent {\bf Proof. }Suppose that $[b_1, b_2]$ contains a point $u$ that belongs to the interior of $\Theta_3$ and $|A(u)|=1$. Then $u\notin\Theta_1$ and $u\notin \Theta_2$. We can choose $v\ne u,\ v\in [b_1, b_2], \ v\notin \Theta_1, \ v\notin \Theta_2$ and $v\in {\rm int} \Theta_3.$ Then $| A(v)|=1.$ We first prove that $u$ is a solution of the problem. Indeed, in this case,
\begin{equation*}
\mathcal{H}(u)=D(u; \Theta_1)+D(u; \Theta_2)+D(u; \Theta_3)=\|b_1-b_2\|-s_1-s_2.
\end{equation*}
Choose $\delta>0$ such that $\B(u; \delta)\cap \Theta_1=\emptyset$ and $\B(u; \delta)\cap \Theta_2=\emptyset$, and $\B(u; \delta)\subseteq \Theta_3$. For every $x\in \B(u; \delta)$, one has
\begin{equation*}
\mathcal{H}(x)=D(x; \Theta_2)+D(x; \Theta_3)=\|x-b_1\|+\|x-b_2\|-s_1-s_2\geq \|b_1-b_2\|-s_1-s_2=\mathcal{H}(u).
\end{equation*}
This implies that $u$ is a local minimum of $\mathcal{H}$, so it is also an absolute minimum since $\mathcal{H}(u)$ is convex. An analogous argument can be applied for $v$. Then $u,\; v\in S.$

Now assume that $S$ has more than one elements. Let $x_*, y_*$ be two distinct elements of $S$. Then $[x_*, y_*]\subseteq S$ by Proposition \ref{pro1}. If there is a solution that does not belong to $\Theta_i$ for every $i\in\{1,2,3\}$, then $S$ reduces to a singleton, that is a contradiction. We can assume without loss of generality that $\Theta_1$ contains infinitely many solutions. Then $\mbox{\rm int }\Theta_1$ contains infinitely many solutions by the strict convexity of $\Theta_1$. If there is such a solution $u$ with $A(u)=\{1\},$ then $u\in\mbox{\rm int }\Theta_1$ and $u$ does not belong to $\Theta_2$, $\Theta_3$. Thus, $u\in [b_2, b_3]$, since if not, the solution must be unique by Lemma \ref{lm2}. Therefore, the conclusion holds. Suppose $|A(u)|=2$ for every solution that belongs to $\mbox{\rm int }\Theta_1$. Then there are infinitely many solutions that lies on the intersection of two sets, which is strictly convex in this case. So there is a solution that belongs to the interior of this intersection, that is a contradiction to Lemma \ref{lm3}. $\h$

\begin{example}{\rm In Figure 1, the given generalized Fermat-Torricelli problem has infinite solutions because $[b_1, b_2]$ intersects the interior of $\Theta_3$, and $|A(u)|=1$ for every $u\in (B,C)$. In fact, the solution set is $S=[B, C]$.
\begin{figure}[!ht]
\centering
\includegraphics[width=7cm]{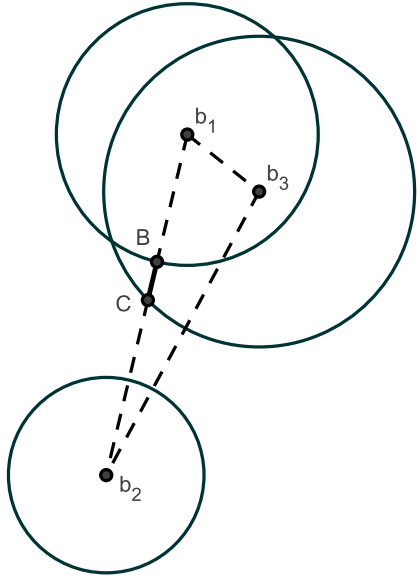}
\caption{A generalized Fermat-Torricelli with infinitely many solutions}\label{fig:}
\end{figure}
}
\end{example}
Let $V_k:=[b_i, b_j]\cap \mbox{\rm int }\Theta_k$ for distinct indices $i,j, k\in \{1,2,3\}.$ As a direct consequence of Theorem \ref{tmu}, we have the following corollary.
\begin{corollary} Suppose $\cap_{i=1}^3\Theta_i=\emptyset$. Problem {\rm (\ref{ft2})} has a unique solution if and only if the following implication holds for any $k\in \{1, 2, 3\}$:
\begin{equation*}
\big[V_k\neq\emptyset\big]\Rightarrow \big[|A(u)|=2 \;\mbox{\rm for all }u\in V_k\big].
\end{equation*}
\end{corollary}
The corollary below provides a sufficient condition for problem (\ref{ft2}) to have a unique solution.
\begin{corollary} Suppose $\cap_{i=1}^3\Theta_i=\emptyset$. Problem {\rm (\ref{ft2})} has a unique solution if $[b_i, b_j]\cap \Theta_k$ contains at most one point for distinct indices $i, j, k\in \{1, 2, 3\}$.
\end{corollary}
\noindent {\bf Proof. }Suppose by contradiction that {\rm (\ref{ft2})} has more than one solutions. By the previous theorem, there is an interval, say $[b_1, b_2]$, that contains a point  belonging to the interior of $\Theta_3$. Then $[b_1, b_2]\cap\Theta_3$ contains infinitely many points, that is a contradiction. $\h$

 \subsection{Solution Constructions}

 In this subsection, we will propose a method of constructing a solution of problem {\rm (\ref{ft2})} for arbitrary balls $\Theta_i$, $i\in J=\{1,2,3\}$.

 \begin{proposition}\label{optim} An element $x_*$ is an optimal solution of problem {\rm (\ref{ft2})} if and only if
 \begin{equation}\label{sc}
 -\sum_{i\in J\setminus A(x_*)}e_i\in \sum_{i\in A(x_*)}[N(x_*; \Theta_i)\cap \B],
 \end{equation}
 where $e_i:=\dfrac{x_*-b_i}{\|x_*-b_i\|}$.
 \end{proposition}
\noindent {\bf Proof. }By the Fermat subdifferential rule (\ref{fermat}), $x_*$ is an optimal solution of problem {\rm (\ref{ft2})} if and only if
\begin{equation*}
0\in \partial \mathcal{H}(x_*)=\sum_{i\in J}\partial D(x_*; \Theta_i)=\sum_{i\in J\setminus A(x_*)}\partial D(x_*;\Theta_i)+\sum_{i\in A(x_*)}\partial D(x_*; \Theta_i).
\end{equation*}
The conclusion then follows from Proposition \ref{subr}. $\h$
\begin{proposition} Consider problem {\rm (\ref{ft2})}. Suppose that $\cap_{i=1}^3\Theta_i\neq \emptyset$. Then
\begin{equation*}
S=\cap_{i=1}^3\Theta_i.
\end{equation*}
\end{proposition}
\noindent {\bf Proof. }Fix any $x_*\in \cap_{i=1}^3\Theta_i\neq \emptyset$. Then $\mathcal{H}(x_*)=0$, and hence $x_*\in S$ since $$\inf_{x\in \Bbb R^2}\mathcal{H}(x)\ge 0.$$ Conversely, for any $x\in S$, one has $\mathcal{H}(x)\le \mathcal H (x_*) =0$. Thus,  $D(x; \Theta_i)=0$ for $i=1,2,3$. This implies $x\in \Theta_i$ for $i=1,2,3$, or equivalently, $x\in \cap_{i=1}^3\Theta_i\neq \emptyset$. $\h$\vspace*{0.05in}

From now on, we only need to consider the case where $\cap_{i=1}^3\Theta_j=\emptyset$. Then $|A(x_*)|< 3$ for all $x_*\in S$.
\begin{proposition} Suppose that $|A(x_*)|=2$ and $x_*$ be an optimal solution of problem {\rm (\ref{ft2})}, say $x_*\in \Theta_1\cap \Theta_2$. Then one of the following condition holds:\\[1ex]
{\rm (1)} $x_*$ is the intersection of $[b_1, b_3]$ and the boundary of $\Theta_1$, and $x_*\in \mbox{\rm int }\Theta_2$.\\
{\rm (2)} $x_*$ is the intersection of $[b_2, b_3]$ and the boundary of $\Theta_2$, and $x_*\in \mbox{\rm int }\Theta_1$.\\
{\rm (3)} $x_*$ belongs to the intersection of the boundary of $\Theta_1$ and the boundary of $\Theta_2$, and \begin{equation*}-e_3=s e_1+te_2,
\end{equation*}
where $s, t\in [0, 1]$ and $e_i=\dfrac{x_*-b_i}{\|x_*-b_i\|}$ for $i=1,2,3$.

Conversely, if one of the conditions above is satisfied, then $x_*$ is an optimal solution of the problem.
\end{proposition}
\noindent {\bf Proof. }Suppose that $|A(x_*)|=2$ and $x_*\in \Theta_1\cap\Theta_2$. By Proposition \ref{optim},
\begin{equation}\label{ot1}
-e_3\in [N(x_*; \Theta_1)\cap \B] + [N(x_*; \Theta_2)\cap \B].
\end{equation}
Since $e_3\neq 0$, one has that $x_*\in \mbox{\rm bd }\Theta_1$ or $x_*\in \mbox{\rm bd }\Theta_2$.\\[1ex]
{\bf Case 1: }$x_*\in \mbox{\rm bd }\Theta_1$ and $x_*\in \mbox{\rm int }\Theta_2$. In this case, (\ref{ot1}) reduces to
\begin{equation*}
-e_3\in N(x_*; \Theta_1).
\end{equation*}
This is equivalent to the fact that $x_*$ is the intersection of $[b_1, b_3]$ and the boundary of $\Theta_1$. Thus, (1) is satisfied.\\[1ex]
{\bf Case 2: }$x_*\in \mbox{\rm bd }\Theta_2$ and $x_*\in \mbox{\rm int }\Theta_1$. Similar to Case 1, under this condition, (2) is satisfied.\\[1ex]
{\bf Case 3: }$x_*\in \mbox{\rm bd }\Theta_1$ and $x_*\in \mbox{\rm bd }\Theta_2$. In this case, (\ref{ot1}) reduces to
\begin{equation*}
-e_3=s e_1+te_2,
\end{equation*}
where $s, t\in [0, 1]$. Notice that in this case, $x_*$ is the intersection of $\mbox{\rm bd }\Theta_1$ and $\mbox{\rm bd }\Theta_2$. Since there are at most two point in this intersection, this condition is verifiable. $\h$
\begin{proposition} Suppose that $|A(x_*)|=1$ and $x_*$ be an optimal solution of problem {\rm (\ref{ft2})}, say $x_*\in \Theta_1$. Then
\begin{equation*}
-e_2-e_3\in N(x_*; \Theta_1) \mbox{ and }\la e_2, e_3\ra \leq -1/2.
\end{equation*}
Conversely, if this condition is satisfied, then $x_*$ is a solution of problem {\rm (\ref{ft2})}.
\end{proposition}
\noindent {\bf Proof. }Condition (\ref{sc}) can be written as
\begin{equation*}
-e_2-e_3\in N(x_*; \Theta_1) \mbox{ and }\|e_2+e_3\|\leq 1.
\end{equation*}
Notice that $\|e_2+e_3\|\leq 1$ if and only if $\la e_2, e_3\ra \leq -1/2$. The proof is complete. $\h$

The construction to find a solution is as follows. We first will find all possible solutions in each set $\Theta_i$ for $i=1,2,3$. If a solution is found, then there is no solution outside of the sets by Lemma \ref{lm1}. If no solution is found, then the optimal solution is found by solving the classical Fermat-Torricelli generated by $b_i$ for $i=1,2,3.$ For instance, solutions on $\Theta_1$ can be found by the following steps:\\[1ex]
{\bf Step 1: }If $\cap_{i=1}^3\Theta_i\neq \emptyset$, then every point in this intersection is a solution. Then, go to the next step.\\
{\bf Step 2: }Connect the centers $[b_1, b_j]$, $j=2,3$. If, for instance, $[b_1, b_3]$ intersects $\mbox{bd }\Theta_1$ at a point that belongs to $\mbox{\rm int }\Theta_2$, then that point is a solution. Then, go to the next step.\\[1ex]
{\bf Step 3:} Find the intersection of the boundary of the balls $\Theta_1$ and $\Theta_2$, and $\Theta_1$ and $\Theta_3$. For instance, let $p$ and $q$ be the intersections of $\mbox{\rm bd }\Theta_1$ and $\mbox{\rm bd }\Theta_2$. Then verify if the condition $-e_3\in [0,1]e_1+[0,1]e_2$ is satisfied at each point. If the condition is satisfied, for instance, at $p$, then $p$ is a solution. Then, go to the next step.\\[1ex]
{\bf Step 4: }Verify if $[b_2, b_3]$ intersects $\Theta_1$. If, for instance, $[b_2, b_3]$ intersects $\Theta_1$, then find all point in the intersection does not belong to $\Theta_2$ and $\Theta_3$, and those are solutions. In the case, $[b_2, b_3]$ does not intersect $\Theta_1,$ then we find $u_2:=[b_1, b_2]\cap \mbox{\rm bd }\Theta_1$, and $u_2:=[b_1, b_3]\cap \mbox{\rm bd }\Theta_1$. Next, find a unique point $q_1$ on the minor curve $u_2u_3$ such that $\widehat{b_2q_1b_1}=\widehat{b_3q_1b_1}$. If $q_1$ does not belong to $\Theta_2$ and $\Theta_3$, and $\widehat{b_2q_1b_3}\geq 120^\circ$, then $q_1$ is a solution of the problem.

The same process can be repeated to find solutions on $\Theta_2$ and $\Theta_3$.

\section{Concluding Remarks}
\label{s:Concluding}

In this paper, we have provided a detailed theoretical analysis for the generalized Sylvester problem and the generalized Fermat-Torricelli problem for Euclidean balls. A natural question is: can we develop numerical algorithms to solve the generalized Sylvester and the generalized Fermat-Torricelli problems for Euclidean balls proposed in this paper? The well-known subgradient method provides such a simple algorithm since subgradients of the distance functions involved in the problems can be explicitly determined. By exploiting the majorization-minimization (MM) principle of computational statistics, another algorithm for solving the constrained version of the generalized Fermat-Torricelli has been developed in \cite{cl}. Faster algorithms can also be developed using other methods of nonsmooth optimization. We will address this in our future research.

\paragraph{Acknowledgments.}  The research of Nguyen Mau Nam was partially
supported by the Simons Foundation under grant \#208785. The research of Nguyen Hoang was partially supported by the NAFOSTED, Vietnam, under grant \# 101.01-2011.26.

\end{document}